\documentclass[11pt]{amsart}
\usepackage{amscd,amssymb}
\usepackage{amsmath}
\usepackage{amsthm,mathtools}
\usepackage{float}
\usepackage{pgf,tikz-cd}
\usepackage{mathrsfs}
\usepackage{fancyhdr}
\usepackage{tabularray}
\usepackage{longtable}
\usepackage{supertabular}
\usepackage{setspace}
\usepackage{color}
\usepackage[matrix,arrow,curve]{xy}
\usepackage{graphicx}
\usepackage{enumerate}
\usepackage{amsfonts}
\usepackage{cite}
\usepackage{pifont}
\usepackage{multirow}
\usepackage{titlesec}
\usepackage{caption}
\usepackage{hyperref}
\usetikzlibrary{arrows}

\newtheorem{theorem}[equation]{Theorem}
\newtheorem*{maintheorem}{Main Theorem}

\newtheorem{lemma}[equation]{Lemma}
\newtheorem{proposition}[equation]{Proposition}
\newtheorem{corollary}[equation]{Corollary}
\newtheorem{conjecture}[equation]{Conjecture}

\theoremstyle{definition}

\newtheorem{definition}[equation]{Definition}

\newtheorem*{ack}{Acknowledgement}

\theoremstyle{remark}

\makeatletter\@addtoreset{equation}{section} \makeatother

\newcommand{\mult}{\operatorname{mult}}
\newcommand{\supp}{\mathrm{Supp}}

\newcommand{\vol}{\operatorname{vol}}

\hypersetup{
  colorlinks   = true,    % Colours links instead of ugly boxes
  urlcolor     = blue,    % Colour for external hyperlinks
  linkcolor    = cyan,    % Colour of internal links
  citecolor    = cyan      % Colour of citations
}

\makeatletter
\def\@seccntformat#1{\@ifundefined{#1@cntformat}%
    {\csname the#1\endcsname\quad}%      default
    {\csname #1@cntformat\endcsname}}%   individual control
\newcommand{\section@cntformat}{\S\thesection~~~}
\makeatother

\title{K-unstable singular del Pezzo surfaces\\ without anticanonical polar cylinder}

\author{In-Kyun Kim}

\address{ \emph{In-Kyun Kim}\newline \textnormal{June E Huh Center for Mathematical Challenges, Korea Institute for Advanced Study, \newline 
\medskip 85 Hoegiro Dongdaemun-gu, Seoul 02455, Republic of Korea.\newline 
\texttt{soulcraw@kias.re.kr}}}%

\author{Jaehyun Kim}

\address{ \emph{Jaehyun Kim}\newline \textnormal{Department of Mathematics, Ewha Womans University \newline 
\medskip 52, Ewhayeodae-gil, Seodaemun-gu, Seoul, 03760, Republic of Korea.\newline 
\texttt{kjh6691@ewha.ac.kr}}}

\author{Joonyeong Won}

\address{ \emph{Joonyeong Won}\newline \textnormal{Department of Mathematics, Ewha Womans University \newline 
\medskip 52, Ewhayeodae-gil, Seodaemun-gu, Seoul, 03760, Republic of Korea.\newline 
\texttt{leonwon@ewha.ac.kr}}}%

\begin{document}

\begin{abstract}

     We prove the existence of singular del Pezzo surfaces that are neither K-semistable nor contain any anticanonical polar cylinder.

\bigskip
\begin{center}
\end{center}

\end{abstract}

%\subjclass[2000]{14R20}

\maketitle

All varieties of this article are assumed to be defined over an algebraically closed field of characteristic zero.

%%%%%%%%%%%%%%%%%%%%%%%%%%%%%%%%%%%%%%%%%%%%%%%%%%%%%%%%%%%%%%%%%%%%%%%%%%%%%%%%%%%%%%%%%

\section{Introduction}
K-stability of Fano varieties is an important research topic that is related to various fields, ranging from K$\ddot{\text{a}}$hler-Einstein metrics to K-moduli spaces. 
In particular, K-stability of birationally superrigid Fano varieties has been extensively studied in \cite{Odaka2013,Kim2018,Stibitz2019,Kim2021b,Kim2023a}. It is conjectured that birationally rigid Fano varieties are always K-stable in \cite{Odaka2013,Kim2018} and the $\alpha$-invariants of those varieties are considered in \cite{Kim2018,Kim2021b,Kim2023a}.
These studies suggest a potential correlation between the rationality of Fano varieties and their K-stability. Considering the connection between rationality and cylindricity, along with the observation that all previously identified non K-semistable Fano varieties are cylindrical, it is natural to speculate the following in \cite{Cheltsov2021}.

\begin{conjecture}[\hspace{1sp}{\cite[Conjecture 1.32]{Cheltsov2021}}]\label{conj}
    Let $X$ be a Fano variety that has at worst klt singularities. If $X$ does not contain $(-K_X)$-polar cylinders, then $X$ is K-polystable.
   \end{conjecture}

A \emph{Cylinder} in a normal projective variety is a Zariski open subset which is isomorphic to $\mathbb{A}^1 \times Z$ for some affine variety $Z$. If the complement of the cylinder is defined by the support of an effective $\mathbb{Q}$-divisor which is $\mathbb{Q}$-linearly equivalent to an ample divisor $H$, then it is said to be $H$-\emph{polar}.
As we can see \cite[Theorem 1.26]{Cheltsov2021}, a cylindrical variety has small $\alpha$-invariant. In this paper, we give an answer to Conjecture 1.1 by considering certain singular del Pezzo surfaces in three-dimensional weighted projective spaces.\\
There are more results that support the conjecture. 
 It has been verified that the automorphism group of a K-polystable Fano variety with klt singularities is reductive in \cite{Alper2020}. This contrasts with the cylindrical view point of unipotent group action.

\begin{proposition}[\hspace{1sp}{\cite[Corollary 3.2]{Kishimoto2013}}]
 For a normal projective variety $X$ with an ample polarization $H$, the affine cone 
 $$\mathrm{Spec}\left(\bigoplus_{m} \mathrm{H}^0(X,\mathcal{O}_{X}(mH))\right)$$
 admits a nontrivial unipotent group action if and only if $X$ has an $H$-polar cylinder.
\end{proposition}

In addition, the conjecture holds at least for smooth del Pezzo surfaces by \cite{Tian1990}. The K-stability of smooth del Pezzo surfaces has been determined as well in \cite{ParkWon2018} by estimating their $\delta$-invariants. 
From now on, we refer to K-unstable if it is not K-semistable.

\begin{proposition}[\hspace{1sp}\cite{{Tian1990}}]
 Let $S$ be a smooth del Pezzo surface of degree $d$. Then the pair $(S,-K_S)$ is
 \begin{itemize}
\item[$\bullet$] $K$-stable, if $d \leq 5$;
  
\item[$\bullet$] $K$-semistable, if $d=6, 9$ or $S \cong \mathbb{P}^1 \times \mathbb{P}^1$;
 
 \item[$\bullet$] $K$-unstable, if $d=7$ or $S \cong \mathbb{F}_1$.
  \end{itemize}

\end{proposition}

\begin{proposition}[\hspace{1sp}{\cite{Cheltsov2016, Cheltsov2017}}]
Let $S$ be a smooth del Pezzo surface of degree $d$. Then $S$ does not contain any $(-K_S)$-polar cylinder if and only if $d \leq 3$.
\end{proposition}  

Meanwhile, quasi-smooth well-formed del Pezzo surfaces in $\mathbb{P}(a_0, a_1,a_2,a_3)$ are thoroughly described by Paemurru in \cite{Paemurru2018}. Let $S$ be a hypersurface in $\mathbb{P}(a_0,a_1,a_2,a_3)$ given by quasi-homogeneous polynomial $\varphi(x,y,z,t)$ of degree $d$ with respect to weights, $1 \leq  a_0 \leq a_1 \leq a_2 \leq a_3$. 
The equation $$\varphi(x,y,z,t) = 0 \subset \mathrm{Spec}(\mathbb{C}[x,y,z,t])$$ defines a three-dimensional hypersurface quasi-homogeneous singularity $(V,o)$, where $o = (0, 0, 0, 0)$. Recall that $S$ is said to be \emph{quasi-smooth} if the singularity $(V, o)$ is isolated and~$S$ is said to be \emph{well-formed} if for each distinct $i, j, k$ and $m, n$ in $\{ 0, 1, 2, 3 \}$,
$$\mathrm{gcd}(a_i, a_j, a_k) = 1~ \hspace{1 cm} ~\mathrm{and}~ \hspace{1 cm} ~\mathrm{gcd}(a_m, a_n) ~\big \vert ~d.$$ 
Moreover, the number defined by $$a_0+a_1+a_2+a_3-d$$
is called the \emph{index} of $S$, and it is known that $S$ is a del Pezzo surface, if $d<a_0+a_1+a_2+a_3$.\\
 The conjecture is anticipated to hold for indices 1 through 6, whereas counter examples are found in very special families of infinite series in \cite{Paemurru2018}. We make a list of K-unstable members in Proposition \ref*{prop:K-unstable} and also prove the absence of cylinders in Theorem \ref*{main2}. For convenience, we present the complete list of infinite series in Table 1 of Section 5. Now we reach the stage to state our main theorem that answer Conjecture \ref{conj}.

\begin{maintheorem}
There are infinitely many quasi-smooth well-formed del Pezzo surfaces with klt singularities in three-dimensional weighted projective spaces which are K-unstable and do not contain any anticanonical polar cylinder.
\end{maintheorem}

The present paper fully understands $(-K_S)$-polar cylindricity and K-stability for the families. We expect that the infinite series exhibit more distinctive properties and it warrants an in-depth investigation.

%%%%%%%%%%%%%%%%%%%%%%%%%%%%%%%%%%%%%%%%%%%%%%%%%%%%%%%%%%%%%%%%%%%%%%%%%%%%%%%%%%%%%%%%%

\section{Preliminaries}
This section is devoted to reminding prerequisite lemmas containing some local inequalities.\\
Let $S$ be a projective surface with at worst klt singularities and let $D$ be an effective $\mathbb{Q}$-divisor on $S$ written by $$D=\sum\limits_{i=1} ^{r} a_i D_i.$$ 
 
\begin{lemma}[\hspace{1sp}{\cite[Proposition 9.5.13]{Lazarsfeld2004a}}]\label{lem:multiplicity}
Let $\mathsf{p}$ be a smooth point of the surface $S$. If the log pair ~$(S,D)$ is not log canonical at $\mathsf{p}$, then $$\mult_{\mathsf{p}}(D) > 1.$$
\end{lemma}

For orbifolds, we have the following analog.

\begin{lemma}[\hspace{1sp}{\cite[Proposition 3.16]{Kollar1997}}]\label{lem:multiplicity_orbifold}
    Let $\mathsf{p}$ be a singular point of type $\frac{1}{r}(a, b)$ of the surface $S$.
    If the log pair $(S,D)$ is not log canonical at $\mathsf{p}$, then $$\mult_{\mathsf{p}}(D) > \frac{1}{r}.$$
    \end{lemma}

\begin{lemma}[\hspace{1sp}{\cite[Theorem 7.5]{Kollar1997}}]\label{lem:adjuction} 
    Let $\mathsf{p}$ be a singular point of type $\frac{1}{r}(a, b)$ of the surface $S$.
    Suppose that the log pair $(S,D)$ is not log canonical at $\mathsf{p}$.\\ If a component ~$D_j$ with $a_j \leq 1$ is smooth at $\mathsf{p}$, then
    $$D_j \cdot (D-a_jD_j)=D_j\cdot\left(\sum_{i \neq j}a_iD_i\right)\geq \sum_{i \neq j}a_i (D_j \cdot D_i)_\mathsf{p} > \frac{1}{r}.$$
    \end{lemma}

\begin{lemma}[\hspace{1sp}{\cite[Lemma 2.2]{Cheltsov2016}}]\label{lem:convex}
Let $T$ be an effective $\mathbb{Q}$-divisor on $S$ such that
\begin{itemize}
 \item[$\bullet$] $T \sim_{\mathbb{Q}} D$ but $T \neq D$;

 \item[$\bullet$] $T=\sum\limits_{i=1}^{r} b_iD_i$ for some non-negative rational numbers $b_i$'s.
\end{itemize}
For every non-negative rational number $\varepsilon$, put $D_{\varepsilon} = (1+\varepsilon)D-\varepsilon T$. Then
\begin{enumerate}
 \item[1.] $D_{\varepsilon} \sim_{\mathbb{Q}} D$ for every $\varepsilon \geq 0$;

 \item[2.] the set $\{ \varepsilon \in \mathbb{Q}_{>0} ~\vert~ D_{\varepsilon} ~\mathrm{is ~effective} \}$ attains a maximum $\mu$;
 
 \item[3.] the support of the divisor $D_{\mu}$ does not contain at least one component of $\supp(T)$;
 
 \item[4.] if $(S,T)$ is log canonical at $\mathsf{p}$ but $(S,D)$ is not log canonical at $\mathsf{p}$, then $(S,D_{\mu})$ is not log canonical at $\mathsf{p}$.
 \end{enumerate}

\end{lemma}

Moreover, if $S$ is rational that contains a $(-K_{S})$-polar cylinder $U \cong \mathbb{A}^1 \times Z$ defined by~$D$, then the natural projection $U \to Z$ induces a rational map $\rho: S \dashrightarrow \mathbb{P}^1$. Let $\mathscr{L}$ be the linear system corresponding to $\rho$.

\begin{lemma}[\hspace{1sp}{\cite[Lemma A.3]{Cheltsov2016}}]\label{lem:boundary of cylinder}
Assume that the base locus of $\mathscr{L}$ is not empty. Then for any effective ~$\mathbb{Q}$-divisor $H \sim_{\mathbb{Q}} -K_{S}$ with $\supp(H) \subset \supp(D)$, the log pair $(S,H)$ is not log canonical at the base point.
\end{lemma}

%%%%%%%%%%%%%%%%%%%%%%%%%%%%%%%%%%%%%%%%%%%%%%%%%%%%%%%%%%%%%%%%%%%%%%%%%%%%%%%%%%%%%%%%%

\section{K-stability}
K-stability for quasi-smooth well-formed del Pezzo surfaces in three-dimensional weighted projective spaces with low indices is completely known by \cite{Johnson2001, Araujo2002, Cheltsov2010, Cheltsov2021a, Kim2022, Kim2023}. Indeed, Johnson and Koll$\acute{\text{a}}$r proved the existence of the K$\ddot{\text{a}}$hler-Einstein metric for almost cases of index $I=1$ using $\alpha$-invariant and Tian's criterion.

\begin{definition}
    Let $X$ be a $\mathbb{Q}$-Fano variety, let $H$ be an ample line bundle on $X$.\\ 
    The $\alpha$-invariant of the log pair $(X,H)$ is the number defined by
    $$ \alpha (X,H) = \mathrm{sup} \left \{ 
        \lambda \in \mathbb{Q} ~\Bigg \vert 
        \begin{array}{ll} 
        &\hspace{-0.3cm} \mathrm{the ~~log ~~pair} ~~(X, \lambda D) \mathrm{~~is ~~log ~~canonical ~~for ~~every} \\
        &\hspace{-0.3cm} \mathrm{effective} ~~\mathbb{Q}\text{-}\mathrm{divisor} ~~D ~~\mathrm{on} ~~X \mathrm{~~with} ~~D\sim_{\mathbb{Q}} H.  
        \end{array}  \right \}. $$
In particular, the invariant $\alpha(X,-K_{X})$ is denoted simply by $\alpha(X)$.   
    \end{definition}

\begin{theorem}[\hspace{1sp}{\cite{DEMAILLY2001,Nadel1990,Tian1987}}]
   Let $X$ be a $\mathbb{Q}$-Fano variety. If $\alpha(X) > \frac{\mathrm{dim}(X)}{\mathrm{dim}(X)+1}$, then $X$ admits an orbifold K$\ddot{a}$hler-Einstein metric.
\end{theorem}

Afterward, the K-stability for remaining ones of $I=1$ is determined in \cite{Araujo2002, Cheltsov2010, Cheltsov2021a}. In particular, the exact values of global log canonical thresholds of infinite series are computed in \cite{Cheltsov2010}. In addition, the K-stability for higher indices, $I=2, 3$, is proved in \cite{Kim2021a, Kim2022, Kim2023} by considering $\delta$-invariants that was introduced by Fujita and Odaka for Fano varieties. For simplicity, we take the Theorem C in \cite{Blum2020} as a definition of $\delta$-invariant.

\begin{definition}
Let $X$ be a $\mathbb{Q}$-Fano variety, let $H$ be an ample line bundle on $X$ and let 
~$\pi$ be a proper birational morphism from a normal variety $Y$ to $X$. Then for
$$A(\mathrm{ord}_E) = 1 + \mathrm{ord}_E(K_{Y/X}),$$
$$S(\mathrm{ord}_E) = \frac{1}{\vol(H)}\int_{0}^{\infty} \vol(\pi^{\ast}H-tE)~dt ,$$
the $\delta$-invariant of the pair $(X,H)$ is the number defined  by
$$\delta(X,H)=\underset{E}{\mathrm{inf}}\dfrac{A(\mathrm{ord}_{E})}{S(\mathrm{ord}_{E})},$$
where $E$ runs through prime divisors over $X$.
The invariant $\delta(X,-K_{X})$ is denoted simply by ~$\delta(X)$.
\end{definition}

The relation of those invariants, $\alpha$ and $\delta$, for big line bundles is proved by Blum in \cite{Blum2020}. Here we present a stronger version that is needed for our results.

\begin{theorem}[\hspace{1sp}{\cite[Theorem A]{Blum2020}}]\label{section 3:alpha and delta inequality}
Let $X$ be a $\mathbb{Q}$-Fano variety. Then for any ample line bundle ~$H$, $\delta(X,H)$ exists and 
\begin{equation*}
    \left(\frac{\mathrm{dim}(X) + 1}{\mathrm{dim}(X)}\right)\alpha(X, H) \leq \delta(X, H)\leq (\mathrm{dim}(X) + 1)\alpha(X, H).
\end{equation*}
Further, the numbers $\alpha(X,H)$ and $\delta(X,H)$ are strictly positive and only depend on the numerical equivalence class of $H$.
\end{theorem}

Through the following result, we can state the criterion by $\delta$-invariant in \cite{Blum2020}.

\begin{theorem}[\hspace{1sp}{\cite[Theorem 1.6]{LiuXuZhuang2022}}]
A log Fano pair $(X, \Delta)$ is uniformly K-stable if and only if it is K-stable.
\end{theorem}

\begin{theorem}[\hspace{1sp}{\cite[Theorem B]{Blum2020}}]
 Let $X$ be a $\mathbb{Q}$-Fano variety.
 \begin{itemize}
  \item[$(i)$] $X$ is K-semistable if and only if $\delta(X) \geq 1$;
 
  \item[$(ii)$] $X$ is K-stable if and only if $\delta(X) > 1$.
  \end{itemize}
\end{theorem}

Along with the above criteria by $\alpha$ or $\delta$ invariants, we can summarize the K-stability of infinite series known so far.  

\begin{lemma}
Let $No.~k$ be the k-th entry of Table 1. Then the followings are hold.    
\begin{itemize}
\item[$(i)$] ($I=1$) $No.~1$ - $6$ with $n=1$ and $No.~23$ are K-stable. 
    
\item[$(ii)$] ($I=2$) 
\begin{itemize}
  \item[$\bullet$] $No.~1$ - $5$ with $n=2$, $No.~8$ with $n=1$ 
 and $No.~24$ - $30$ are K-stable;

 \item[$\bullet$] $No.~6$ with $n=2$ and  $No.~7$ with $n=1$ are K-unstable.
 \end{itemize}

\item[$(iii)$] ($I=3$) $No.~1$ - $6$ with $n=3$ is K-unstable. 

\item[$(iv)$] ($I=4,6$) $No.~31$ - $35$ are K-stable. 
\end{itemize}
\end{lemma}

\begin{proof}
For $(i)$, refer to \cite[Theorem 8]{Johnson2001}, \cite[Theorem 4.1]{Araujo2002}, \cite[Theorem 1.10]{Cheltsov2010} and \cite[Theorem 1.7]{Cheltsov2021a}. In \cite{Cheltsov2010}, it is also proved that the log canonical thresholds of $No.~31$ - $35$ are $1$ which implies $(iv)$. In addition, for indices $I=2, 3$, see \cite{Boyer2003,Cheltsov2010, Cheltsov2013a} and the $\delta$-invariants are computed in \cite{Kim2022, Kim2023}.
\end{proof}

There are more K-unstable members for sufficiently large indices. 

\begin{proposition}\label{prop:K-unstable}
Let $No.~k$ be as above. Then all of the followings are K-unstable.
\begin{itemize}
 \item[$\bullet$] $No.~1$ - $6$, $No.~8$ and $No.~11$ with $n>3$;

 \item[$\bullet$] $No.~7$ with $n>1$;
 
 \item[$\bullet$] $No.~9,~10$ and $No.~12$ - $22$ with $n>2$. 
\end{itemize}

\end{proposition}

\begin{proof}
Let $S_k$ be a quasi-smooth member of the infinite series $No.~k$ with the weight of $x$ is $a_0$ and let $I$ be the index of $S_k$. If $\frac{a_0}{I} < \frac{1}{3}$, then by Theorem \ref{section 3:alpha and delta inequality} we have the following.
$$\delta(S_k) \leq 3 \alpha(S_k) \leq  \dfrac{3a_0}{I} < 1.$$
These imply that $S_k$ is K-unstable
\end{proof}

K-stability of infinite series is presented in Table 2 of Section 5.

%%%%%%%%%%%%%%%%%%%%%%%%%%%%%%%%%%%%%%%%%%%%%%%%%%%%%%%%%%%%%%%%%%%%%%%%%%%%%%%%%%%%%%%%%

\section{Absence of Cylinders}
In this section, we prove that none of quasi-smooth members from thirty-five families of infinite series contain anticanonical polar cylinders. From the following result and \cite{Cheltsov2010} we find K-polystable Fano varieties which do not contain $(-K_X)$-polar cylinders.

%Before delving into our main assertion, we reduce the cases by \cite{Kishimoto2014} and \cite{Cheltsov2010}.

\begin{theorem}[\hspace{1sp}{\cite[Theorem 1.26]{Cheltsov2021}}]\label{thm:alpha obstruction}
    Let $X$ be a Fano variety that has at most klt singularities. If $\alpha(X)\geq 1$, then $X$ does not contain $(-K_X)$-polar cylinder.
\end{theorem}

% \begin{proposition}[\hspace{1sp}{\cite[Lemma 4.8]{Kishimoto2014}}]
% Let $S$ be a smooth del Pezzo surface and let $D \sim_{\mathbb{Q}} -K_{S}$ be an effective $\mathbb{Q}$-divisor on $S$. If $S$ admits a $(-K_{S})$-polar cylinder by $D$, then $(S,D)$ is not log canonical at some point $\mathsf{P}$.
% \end{proposition}

% As we will see in the proof of Lemma 4.7, this holds for our singular del Pezzo surfaces.

% \begin{corollary}
% If $S$ admits $(-K_S)$-polar cylinder, then $\alpha(S) < 1$.
% \end{corollary}

The $\alpha$-invariants of quasi-smooth members of infinite series are computed in \cite{Cheltsov2010}. In particular, the $\alpha$-invariants of quasi-smooth members of $No.~23$ - $35$ are all one. Based on this, we can readily conclude that they do not have any anticanonical polar cylinder. Therefore, for the remainder of this paper, we only need to consider quasi-smooth members of $No.~1$ - $22$.

\subsection{Non-lc locus of a pair}
As a first step, we find a section of the surface that plays the same role of in \cite[Theorem 1.12]{Cheltsov2016} without any assumption on the existence of cylinder. We provide details for only some entries of the list, while those of all other cases are recorded in Section 5. The contents corresponding to the proof of Lemma \ref*{Lemma:lc for S_k-H_x} (resp. Theorem \ref*{main1}) are listed in Table 3 (resp. Tables 5, 6, 7). We begin with summarizing our notations:

\begin{itemize}
\item[$\bullet$] Let $\mathbb{P}(a_0, a_1, a_2, a_3)_{x,y,z,t}$ be the weighted projective space with variables $x, y, z, t$ of weights $a_0, a_1,a_2,a_3$, respectively.

\item[$\bullet$] Let $S_k$ be the quasi-smooth well-formed del Pezzo surface in $\mathbb{P}(a_0, a_1, a_2, a_3)_{x,y,z,y}$ which is listed in $k$-th entry of Table 1.

\item[$\bullet$] Let $D\sim_{\mathbb{Q}} -K_{S_k}$ be an effective $\mathbb{Q}$-divisor on $S_k$.

\item[$\bullet$] $H_x$ and $H_y$ are the hyperplane sections of $S_k$ defined by $x=0$ and $y=0$, respectively.

\item[$\bullet$] $\mathsf{p}_{_x} = [1:0:0:0], ~\mathsf{p}_{_y}=[0:1:0:0], ~\mathsf{p}_{_z}=[0:0:1:0],$ and $\mathsf{p}_{_t}=[0:0:0:1]$. 

\item[$\bullet$] For a positive integer $m$, $\mathsf{p}_{_m}$ is a cyclic quotient singular point other than $\mathsf{p}_{_x}$ through $\mathsf{p}_{_t}$ which is induced by the primitive $m$-th root of unity.

\end{itemize}

For the proof of Lemma \ref*{Lemma:lc for S_k-H_x}, it is worthwhile to present one result in \cite{Araujo2002}.

\begin{lemma}[\hspace{1sp}{\cite[Corollary 3.7]{Araujo2002}}]\label{lem:supportM}
    Let $S$ be an anticanonically embedded quasi-smooth log del Pezzo surface of degree $d$ in $\mathbb{P}(a_0, a_1, a_2, a_3)_{x,y,z,t}$. Let $\pi_t : X \to \mathbb{P}(a_0, a_1, a_2)$ denote the projection from $\mathsf{p}_{_t} = (0: 0: 0: 1)$. Assume that $\pi_t$ has only finite fibers. If $H^0(\mathbb{P}, \mathcal{O}_{\mathbb{P}}(r))$ contains

\begin{itemize}
    \item[$\bullet$] at least two different monomials of the form $x^{i}y^{j}$,
    \item[$\bullet$] at least two different monomials of the form $x^{\ell}z^{m}$, 
\end{itemize}
    
    where $r$ is a positive integer and the four constants $i,j, \ell$ and $m$ are non-negative integers, then for every $\mathsf{p} \in S \setminus (x = 0)$ and every effective $\mathbb{Q}$-divisor $D \sim_{\mathbb{Q}} -K_S$,
    there is a divisor $F \in \vert \mathcal{O}_{S}(r) \vert$ that can be written as $F = F_0 + a(x = 0)$, where $a$ is a nonnegative integer, and $F_0$ is an effective Weil divisor passing through $\mathsf{p}$ and not containing any irreducible component of $D$. Hence, 
    $$ \dfrac{r d}{a_0 a_1 a_2 a_3} \geq \mathrm{mult}_{\mathsf{p}}(D).$$
\end{lemma}

% \begin{lemma} For each positive integer $k \leq 22$, if $(S_k,D)$ is not log canonical at $\mathsf{P}$, then $\mathsf{P}$ is contained in the section $H_x$ defined by $x=0$.
% \end{lemma}
Let $D\sim_{\mathbb{Q}} -K_{S_k}$ be an effective $\mathbb{Q}$-Cartier divisor such that the log pair $(S_k, D)$ is not log canonical.
\begin{lemma}\label{Lemma:lc for S_k-H_x}
    For each positive integer $k \leq 22$, the log pair $(S_k, D)$ with $n>2$ is log canonical along $S_k\setminus H_x$.
\end{lemma}

\begin{proof} 
% Suppose on the contrary, $\mathsf{P}$ is not contained in $\supp(H_x)$.
%Suppose that the log pair $(S_k, D)$ is not log canonical at a point $\mathsf{P}\in S_k\setminus H_x$.
We provide the details only for $S_1$ and $S_{22}$, since the proof of all other cases are almost verbatim.

\begin{itemize}
    \item[$\bullet$] $S_{1} \subset \mathbb{P}(1,3n-2,4n-3,6n-5)$ of degree $12n-9$ with index $n$;

    \item[$\bullet$] $S_{22} \subset \mathbb{P}(7,28n+6,42n+9,63n+10)$ of degree $126n+27$ with index $7n+5$.
\end{itemize}

By a suitable coordinate change, we may assume that $S_1$ (resp. $S_{22})$ are defined by a quasi-homogeneous polynomial
$$x^{12n-9}+z^3+xt^2+y^2t=0 ~~(\mathrm{resp.} ~~xt^2+z^3+y^3z+x^{14n+3}y =0).$$

\textbf{Case 1.} Suppose that the log pair $(S_1,D)$ is not log canonical at $\mathsf{p} \in S_1 \setminus H_x$. Then $S_1$ is smooth at $\mathsf{p}$. Let $\mathscr{M}$ be the linear system defined by the monomials
\begin{equation*}
    x^{4n-3}, z, x^{n-1}y
\end{equation*}
Then there is a member $M \in \mathscr{M}$ such that $\mathsf{p} \in M$ and $M \nsubseteq \supp(D)$ by Lemma \ref*{lem:supportM}. We have
$$\dfrac{n(12n-9)}{(3n-2)(6n-5)} = D\cdot M \geq \mult_{\mathsf{p}} (D) > 1.$$
It is impossible for $n>2$. 

\textbf{Case 2.} Suppose that the log pair $(S_{22},D)$ is not log canonical at $\mathsf{p} \in S_{22}\setminus H_x$. Then $\mathsf{p}$ is either a smooth point on $S_{22}$ or $\mathsf{p}=\mathsf{p}_{_x}$. Let $\mathscr{M}$ be the linear system defined by the monomials
\begin{equation*}
    x^{28n+6}, y^7
\end{equation*}
% $$\Big \vert ~\alpha x^{28n+6}+\beta y^7 =0~ \Big \vert,$$
Then there is a member $M \in \mathscr{M}$ such that $\mathsf{p} \in M$. 
% Since $\mathsf{P} \not\in H_x$, there are two subcases as follows:
% \begin{itemize}
%  \item[$\bullet$] $\alpha=0$, $M=H_y$;

%  \item[$\bullet$] $\alpha \neq 0, \beta \neq 0$.
%  \end{itemize}
We first consider the case that $\mathsf{p}\in H_y$, where $H_y$ is the hyperplane section defined by $y=0$. The hyperplane section $H_y$ is isomorphic to the hypersurface that is embedded in the weighted projective space $\mathbb{P}(7,42n+9,63n+10)_{x,z,t}$ and is given by the quasi-homogeneous polynomial 
$$xt^2+z^3=0.$$
From this, we can see that $H_y$ is irreducible and the log canonical threshold of $(S_{22},H_y)$ is $\frac{5}{6}$.
Hence, the log pair $(S_{22},\frac{7n+5}{28n+6}H_y)$ is log canonical at $\mathsf{p}$. By Lemma \ref*{lem:convex}, we may assume that $$H_y \not\subset \supp(D).$$ This implies a contradiction when $n>2$,
$$\dfrac{21n+15}{441n+70}=H_y\cdot D \geq \mult_{\mathsf{p}}(H_y)\cdot\mult_{\mathsf{p}}(D) >
\left\{ \begin{array}{ll}
1, ~~\mathrm{if} ~\mathsf{p} ~~\mathrm{is ~~smooth},\\\\   
\dfrac{1}{7}, ~~\mathrm{if} ~\mathsf{p}=\mathsf{p}_{_x}.\\
\end{array} \right.$$

We next consider the case that $\mathsf{p} \in M$ and $\mathsf{p} \not\in H_x \cup H_y$. In this situation there are nonzero constants $\alpha$ and $\beta$ such that $M$ is defined by $\alpha x^{28n+6}+\beta y^7 =0$. Then $\mathsf{p}$ is a smooth point of $S_{22}$ and since ~$M$ is given by the quasi-homogeneous polynomial
$$xt^2+z^3+y^3 z+x^{14n+3} y=\alpha x^{28n+6}+\beta y^7 =0,$$
$M$ is irreducible such that $\mult_{\mathsf{p}}(M) \leq 2$. Hence, the log pair $(S_{22},\frac{7n+5}{196n+42}M)$ is log canonical at ~$\mathsf{p}$. We have
$$\dfrac{21n+15}{63n+10}=M\cdot D \geq \mult_{\mathsf{p}}(M)\cdot \mult_{\mathsf{p}}(D) > 1.$$
This also implies a contradiction when $n>2$.
\end{proof}

The singular points and linear systems needed for the proof of other cases are listed in Table 3.

We reach the stage to state our first main result.

\begin{theorem}{\label{main1}}
    Under the same assumption in Lemma \ref{Lemma:lc for S_k-H_x}, the hyperplane section $H_x$ is contained in the support of $D$.
\end{theorem}
    
\begin{proof}
We divide the proof into two cases depending on the irreducibility of the hyperplane section $H_x$. We set 
\begin{equation*}
    \mathrm{IR} = \{1, 3, 4, 5, 6, 8, 9, 10, 12, 13, 15, 17, 19, 21\},\qquad
    \mathrm{RE} = \{2,7,11,14,16,18,20,22\}.
\end{equation*}
The $H_x$ is irreducible for quasi-smooth members of the infinite series $No.~i$, where $i\in \mathrm{IR}$ and is reducible for quasi-smooth members of the infinite series $No.~j$, where $j\in \mathrm{RE}$. We provide details only for $No.~1$ for the irreducible cases and $No.~2$ and  $No.~22$ for the reducible cases. Suppose that $H_x$ is not contained in $\supp(D)$.\\

\textbf{Case 1.} We consider a quasi-smooth member $S_1$ of the infinite series $No.~1$. We can see that the hyperplane section $H_x$ is irreducible. Suppose that $H_x \not\subset \supp(D)$.
We obtain the following inequality
$$\dfrac{3n}{(3n-2)(6n-5)}=H_x\cdot D \geq \mult_{\mathsf{p}}(H_x)\cdot \mult_{\mathsf{p}}(D)> \left\{
        \begin{array}{ll}
        1, ~~\mathrm{if} ~\mathsf{p} ~~\mathrm{is ~~smooth},\\\\   
        \dfrac{2}{3n-2}, ~~\mathrm{if} ~\mathsf{p}=\mathsf{p}_{_y},\\\\
        \dfrac{2}{6n-5}, ~~\mathrm{if} ~\mathsf{p}=\mathsf{p}_{_t},\\
        \end{array} \right.$$
which implies a contradiction when $n > 2$. All other absurd inequalities where $H_x$ is irreducible are listed in Table 5.\\

\textbf{Case 2.} We consider the case that $H_x$ is reducible.
The details are provided only for $S_2$ and $S_{22}$, since the proof of $S_2$ is slightly different from that of $S_j$, where $j\in \mathrm{RE}\setminus \{2\}$.

\begin{itemize}
    \item[$\bullet$] $S_{2} \subset \mathbb{P}(1,3n-2,4n-3,6n-4)$ of degree $12n-8$ with index $n$;
   
    \item[$\bullet$] $S_{22} \subset \mathbb{P}(7,28,+6,42n+9,63n+10)$ of degree $126n+27$ with index $7n+5$.\\
    \end{itemize}

\textbf{Subcase 2.1.} We consider $S_{2} \subset \mathbb{P}(1,3n-2,4n-3,6n-4)$.
By a suitable coordinate change, we may assume that $S_2$ is defined by the quasi-homogeneous polynomial
$$t(t+y^2)-x(z^3+x^{12n-9})=0$$ 
and $H_x$ is written by $$H_x=L_1+L_2,$$ 
where $L_1$ and $L_2$ is defined by $x=t=0$ and $x=t+y^2=0$, respectively. If both of $L_1$ and $L_2$ are not contained in $\supp(D)$, then we obtain the following inequality
$$\dfrac{2n}{(3n-2)(4n-3)}=H_x\cdot D \geq \mult_{\mathsf{p}}(H_x)\cdot \mult_{\mathsf{p}}(D)> \left\{
    \begin{array}{ll}
    1, ~~\mathrm{if} ~\mathsf{p} ~~\mathrm{is ~~smooth},\\\\   
    \dfrac{2}{4n-3}, ~~\mathrm{if} ~\mathsf{p}=\mathsf{p}_{_z},\\\\
    \dfrac{2}{3n-2}, ~~\mathrm{if} ~\mathsf{p}=\mathsf{p}_{_{3n-2}},\\
    \end{array} \right.$$
which implies a contradiction when $n > 2$.\\

If one of $L_1$ and $L_2$ is not contained in $\supp(D)$, without loss of generality we can assume that 
$$L_1 \subset \supp(D) ~~\mathrm{and} ~~L_2 \not\subset \supp(D).$$  
Let $Y\subset \mathbb{P}(1,3n-2,4n-3,6n-4,6n-5)_{x,y,z,t,w}$ be the complete intersection defined by $$wx-t-y^2=wt+z^3+x^{12n-9}=0,$$
and let $\pi : S_2 \to Y$ be the morphism defined by
\begin{equation*}
    (x:y:z:t) \mapsto \left(x:y:z:t:w=\frac{t+y^2}{x}=-\frac{z^3+x^{12n-9}}{t}\right).
\end{equation*}     
Then $Y$ is isomorphic to $S_1$. And $L_1$ is contracted to the point $\mathsf{q}=(0:0:0:0:1)$ by $\pi$. This implies that for an effective $\mathbb{Q}$-divisor $\pi_{\ast}(D)\sim_{\,\mathbb{Q}}-K_{Y}$ and the curve $\pi_{\ast}(L_2)$ defined by $x=0$ on ~$Y$,  the pair $(Y,\pi_{\ast}(D))$ is not log canonical at $\mathsf{q}$ with
$$ \pi_{\ast}(L_2) \not\subset \supp(\pi_{\ast}(D)),$$
which contradicts to Case 1.\\

\textbf{Subcase 2.2.} We consider $S_{22} \subset \mathbb{P}(7,28n+6,42n+9,63n+10)$. By a suitable coordinate change, we may assume that $S_{22}$ is defined by a quasi-homogeneous polynomial
$$xt^2+z^3+y^3z+x^{14n+3}y=0$$ 
and $H_x$ is written by $$H_x=L+R,$$ 
where $L$ and $R$ is defined by $x=z=0$ and $x=z^2+y^3=0$, respectively. We have the following intersection numbers:

$$\emph{(A0)} \hspace{0.5 cm} \small{\begin{array}{ll}
    L\cdot(&\hspace{-0.3 cm}-K_{S_{22}})=\dfrac{7n+5}{(28n+6)(63n+10)}, ~~R\cdot(-K_{S_{22}})=\dfrac{14n+10}{(28n+6)(63n+10)}, ~~L\cdot R=\dfrac{3}{63n+10},\\\\
    &L^2=\dfrac{7}{(28n+6)(63n+10)}-\dfrac{3}{63n+10}, ~~R^2=\dfrac{14}{(28n+6)(63n+10)}-\dfrac{3}{63n+10}.\\\\
   \end{array}}$$

If both of $L$ and $R$ are not contained in $\supp(D)$, then we obtain the following inequalities
$$\emph{(A1)} \small{\hspace{0.5 cm} \dfrac{3(7n+5)}{(28n+6)(63n+10)}=H_x\cdot D \geq \mult_{\mathsf{p}}(H_x)\cdot \mult_{\mathsf{p}}(D)> \left\{
    \begin{array}{ll}
    1, ~~\mathrm{if} ~\mathsf{p} ~~\mathrm{is ~~smooth},\\\\   
    \dfrac{2}{28n+6}, ~~\mathrm{if} ~\mathsf{p}=\mathsf{p}_{_y},\\\\
    \dfrac{2}{63n+10}, ~~\mathrm{if} ~\mathsf{p}=\mathsf{p}_{_t},\\\\
    \dfrac{2}{14n+3}, ~~\mathrm{if} ~\mathsf{p}=\mathsf{p}_{_{14n+3}}.\\
    \end{array} \right.}$$
From these we obtain a contradiction.\\

Now suppose that one of $L$ and $R$ is not contained in $\supp(D)$. If ~$$R\subset \supp(D) ~~(\mathrm{resp}. ~~L \subset \supp(D)),$$ then $D$ is written by 
$$D=\lambda_{R}R+\Delta ~~\mathrm{with} ~~R\not\subset \supp(\Delta)~~(\mathrm{resp.} ~~D=\lambda_{L}L + \Delta  ~~\mathrm{with} ~~L\not\subset \supp(\Delta)).$$
First, for the intersection point $\mathsf{p}_{_t}$ of $L$ and ~$R$, the pair $(S_{22},D)$ is log canonical at $\mathsf{p}_{_t}$ from the following absurd inequalities when $\mathsf{p}=\mathsf{p}_{_t}$
$$\emph{(A2)} \begin{array}{ll}
    &\dfrac{7n+5}{(28n+6)(63n+10)}=L\cdot D > \dfrac{1}{63n+10}\\\\
    & \left(\mathrm{resp}. ~~\dfrac{14n+10}{(28n+6)(63n+10)}=R\cdot D > \dfrac{1}{63n+10}
 \right).\\ 
\end{array}$$

 From this if $R \subset \supp(D) ~~\mathrm{and} ~~L \not\subset \supp(D)$, then $\lambda_R \leq 1$ so that by the inversion of adjunction, we obtain the following inequality
$$\emph{(A3)} \hspace{0.5 cm} \dfrac{14n+10+\lambda_{R}(84n+4)}{(28n+6)(63n+10)}=R\cdot (D-\lambda_{R}R)= R \cdot\Delta > \left\{ 
    \begin{array}{ll}
    1, ~~\mathrm{if} ~\mathsf{p} ~~\mathrm{is ~~smooth},\\\\   
    \dfrac{1}{14n+3}, ~~\mathrm{if} ~\mathsf{p}=\mathsf{p}_{_{14n+3}},\\
    \end{array} \right. $$
which implies a contradiction.\\

Lastly, if $L \subset \supp(D) ~~\mathrm{and} ~~R \not\subset \supp(D)$, then there is a morphism $\pi : Y \to S_{22}$ defined by
\begin{equation*}
    (x:y:z:t) \mapsto \left(x:y:z:t:w=\frac{z^2+y^3}{x}=-\frac{t^2+x^{14n+2}y}{z}\right),
\end{equation*}
where $Y\subset \mathbb{P}(7,28n+6,42n+9,63n+10,84n+11)_{x,y,z,t,w}$ is the complete intersection defined by $$wx-z^2-y^3=wz+t^2+x^{14n+2}y=0$$
with index $7n+5$. Then $L$ is contracted to the point $\mathsf{q}=(0:0:0:0:1)$ by $\pi$ and for an effective $\mathbb{Q}$-divisor $D'=\pi_{\ast}(D)\sim_{\mathbb{Q}}-K_{Y}$,
$$R'= \pi_{\ast}(R) \not\subset \supp(D').$$ From this we obtain the following inequality
$$\emph{(A4)} \hspace{0.5 cm} \dfrac{14n+10}{(14n+3)(84n+11)}=R'\cdot D' \geq \mult_{\mathsf{q}}(R')\cdot \mult_{\mathsf{q}}(D')> \dfrac{3}{84n+11},$$
which also implies a contradiction.\\

All other absurd inequalities corresponding to $\emph{(A0)}$ through $\emph{(A4)}$ where $H_x$ is reducible are listed in Tables 6 and 7.\\
\end{proof}/Users/jaehyunkim/Documents/Screenshot/infinite series 1.png /Users/jaehyunkim/Documents/Screenshot/infinite series 2.png

\subsection{Boundary of Cylinder}
From now on, we additionally assume that $S_k$ admits a $(-K_{S_k})$-polar cylinder defined by an effective $\mathbb{Q}$-Cartier divisor $D$. It means that the Zariski open subset $$S_k \setminus \supp(D)$$ is isomorphic to $\mathbb{A}^1 \times Z$ for some affine curve $Z$. It is also said that $D$ is a boundary of the cylinder. Then we have a following diagram:

\begin{equation}\label{eq:pencil}
\xymatrix{\mathbb{A}^1\times Z\cong U~\ar[dd]^{p_2}\ar@{^{(}->}[rr]& &S_k \ar@{-->}^{\rho}[dd]& &W\ar[ll]_{\psi} \ar[dd]^{\sigma} \ar[ddll]^{\widetilde{\rho}}&  \\ 
& & & \\
Z~\ar@{^{(}->}[rr]& &\mathbb{P}^1 & & \mathbb{F}_1 \ar[ll]&  }
\end{equation}
 
where $p_2$ is the projection to second factor, $\rho$ is the induced rational map with the corresponding pencil $\mathscr{L}$ and $\tilde{\rho}$ is a morphism obtained by resolving the indeterminacy of $\rho$, if any.

\begin{lemma}\label{lem:baseptfree}
The base locus of the pencil $\mathscr{L}$ cannot be empty. 
\end{lemma}

\begin{proof}
Suppose that the linear system $\mathscr{L}$ induced by the cylinder with boundary $$D=\sum_{i=1}^{r} \lambda_i C_i$$ 
defines a conic bundle $\rho$ with the section $C_1$. Then for a generic fiber $F$ of the pencil,
$$-2+\lambda_1 =\left(K_{S_k}+\sum_{i=1}^{r} \lambda_i C_i \right)\cdot F = (K_{S_k}+D)\cdot F = 0.$$
This implies that $(S_k,D)$ is not log canonical along $C_1$ and by Lemma \ref*{Lemma:lc for S_k-H_x}, $$C_1 \subset H_x.$$ 
As before, we provide the details only for $No.~1$ and $No.~22$ depending on the irreducibility of the section $H_x$. For $S_1$, the section $H_x$ is irreducible and is not contained in the support of an effective $\mathbb{Q}$-divisor $$D^{\ast}=\dfrac{n}{n-2}D-\dfrac{2n}{n-2}H_x = \sum_{i=2}^{r} \lambda^{\ast}_i C_i \sim_{\mathbb{Q}} -K_{S_1},$$
where $\lambda^{\ast}_i = \dfrac{n \lambda_i}{n-2}$. This implies a contradiction,
$$-2=\left(K_{S_1}+\sum_{i=2}^{r} \lambda^{\ast}_i C_i\right)\cdot F = (K_{S_1}+D^{\ast})\cdot F = 0.$$

Now for $S_{22}$, since the section $H_x=L+R$ is reducible, $D$ is written by $$D=\lambda_{L}L+ \lambda_{R}R+\Delta,$$
where both of $L$ and $R$ are not contained in $\supp(\Delta)$.
If $C_1=L$ (resp. $C_1=R$), then $$\lambda_L = \lambda_1=2 ~~(\mathrm{resp}. ~~\lambda_R =\lambda_1=2).$$ Hence, $L$ (resp. $R$) is not contained in the support of an effective $\mathbb{Q}$-divisor 
$$D^{\ast}=\dfrac{(7n+5)}{(7n-9)}D-\dfrac{2(7n+5)}{(7n-9)}H_x= \sum_{i=2}^{r} \lambda^{\ast}_i C_i \sim_{\mathbb{Q}} -K_{S_{22}}.$$
This also implies a contradiction,
$$-2=\left(K_{S_{22}}+\sum_{i=2}^{r} \lambda^{\ast}_i C_i\right)\cdot F = (K_{S_{22}}+D^{\ast})\cdot F = 0.$$

\end{proof}
 
All other $D^{\ast}$ are summarized in Table 4. Now we reach the stage to state our second main result.

\begin{theorem}{\label{main2}}
    All thirty-five families of infinite series do not contain any anticanonical polar cylinder, if $n > 2$.
\end{theorem}

\begin{proof}
For a positive integer $k \leq 22$, suppose that $S_k$ admits a $(-K_{S_k})$-polar cylinder defined by an effective $\mathbb{Q}$-divisor $$D=\sum_{i=1}^{r} \lambda_i C_i.$$
Let $\mathsf{p}$ be the base point of the corresponding linear system. Then for a resolution of indeterminacy $\tilde{\rho} : W \to \mathbb{P}^1$ at $\mathsf{p}$, let
$$K_{W}+\sum_{i=1}^{r} \lambda_i \widetilde{C_i}+\sum_{j=1}^{m} \mu_j E_j=\tilde{\rho}^{\ast}(K_{S_k}+D),$$
where $E_1$ is the section in $W$ with $(E_1)^2=-1$. Then for a birational transform $\widetilde{F}$ of $F$ in $W$,
$$-2+\mu_1=\left( K_{W}+\sum_{i=1}^{r} \lambda_i \widetilde{C_i}+\sum_{j=1}^{m} \mu_j E_j\right)\cdot \widetilde{F}=\tilde{\rho}^{\ast}(K_{S_k}+D)\cdot\widetilde{F}=0.$$

This implies that the log pair $(S_k,D)$ is not log canonical at $\mathsf{p}$. Meanwhile, by Lemmas \ref*{lem:convex} and \ref*{lem:boundary of cylinder}, we have an effective $\mathbb{Q}$-divisor $D' \sim_{\mathbb{Q}}-K_{S_k}$ such that the log pair $(S_k,D')$ is not log canonical at $\mathsf{p}$ and 
$$H_x \not\subset \supp(D')$$
which contradicts to Theorem \ref*{main1}. In addition, it has been proved that the surfaces $S_{23}$ through $S_{35}$ do not contain any anticanonical polar cylinder by Theorem \ref*{thm:alpha obstruction}.

\end{proof}

Through the above with the results on K-stability in Section 2, we obtain the following. 

\begin{corollary}
There are infinitely many quasi-smooth well-formed del Pezzo
surfaces with klt singularities in three-dimensional weighted projective spaces which are K-stable and do not contain any anticanonical polar cylinder.
\end{corollary}

\begin{corollary}
There are infinitely many quasi-smooth well-formed del Pezzo
surfaces with klt singularities in three-dimensional weighted projective spaces which are K-unstable and do not contain any anticanonical polar cylinder.
\end{corollary}

%%%%%%%%%%%%%%%%%%%%%%%%%%%%%%%%%%%%%%%%%%%%%%%%%%%%%%%%%%%%%%%%%%%%%%%%%%%%%%%%%%%%%%%%%

\section{Appendix}
\label{sec:5}

Complete list of infinite series in \cite{Paemurru2018} are presented in Table 1. For each positive integer $n$, let~$d$ be the degree of each surface and $I$ be the index of the surface defined by $a_0+a_1+a_2+a_3-d$.

\begin{center}
    \begin{longtable}{ccccccccccccccccccccccc}
    \caption{Infinite Series of del Pezzo Surfaces in $\mathbb{P}(a_0,a_1,a_2,a_3)$}\\
    \hline\\[-1.5 ex]
    \footnotesize{No.} && & &&& & \footnotesize{$(a_0,a_1,a_2,a_3)$}& &&& &&& \footnotesize{$d$}& & &&& && \footnotesize{$I$}& \\ [1 ex]
    \hline\hline\\[-1.5 ex]
    \footnotesize{\emph{1}} & &&&& & & \footnotesize{$(1,3n-2,4n-3,6n-5)$} && &&& & & \footnotesize{$12n-9$}& && &&& & \footnotesize{$n$} &  \\[0.5 ex]
    \hline\\[-1.5 ex]

    \footnotesize{\emph{2}} & &&&& & & \footnotesize{$(1,3n-2,4n-3,6n-4)$} & && && & &\footnotesize{$12n-8$}  & & & &&&& \footnotesize{$n$} & \\[0.5 ex]
    \hline\\[-1 ex]

    \footnotesize{\emph{3}} & & &&&& & \footnotesize{$(1,4n-3,6n-5,9n-7)$} & && && && \footnotesize{$18n-14$} & & &&&& & \footnotesize{$n$} & \\[0.5 ex]
    \hline\\[-1 ex]

    \footnotesize{\emph{4}} & & &&& && \footnotesize{$(1,6n-5,10n-8,15n-12)$} & & &&&& & \footnotesize{$30n-24$} & & &&& && \footnotesize{$n$} & \\[0.5 ex]
    \hline\\[-1 ex]

    \footnotesize{\emph{5}} & & &&& && \footnotesize{$(1,6n-4,10n-7,15n-10)$} & & &&& & &\footnotesize{$30n-20$} & & &&&& & \footnotesize{$n$} & \\[0.5 ex]
    \hline\\[-1 ex]

    \footnotesize{\emph{6}} & & && &&& \footnotesize{$(1,6n-3,10n-5,15n-8)$} & & &&& & &\footnotesize{$30n-15$} & & & &&&& \footnotesize{$n$} &  \\[0.5 ex]
    \hline\\[-1 ex]

    \footnotesize{\emph{7}} & & & &&& & \footnotesize{$(1,8n-2,12n-3,18n-5)$} & & &&& & &\footnotesize{$36n-9$} & & &&&& & \footnotesize{$2n$} &  \\[0.5 ex]
    \hline\\[-1 ex]

    \footnotesize{\emph{8}} & & &&&& & \footnotesize{$(2,6n-3,8n-4,12n-7)$} & & &&&& & \footnotesize{$24n-12$} & & &&&& & \footnotesize{$2n$} &  \\[0.5 ex]
    \hline\\[-1 ex]

    \footnotesize{\emph{9}} & & &&&& & \footnotesize{$(2,6n+1,8n+2,12n+3)$} & &&& && & \footnotesize{$24n+6$} & & && &&& \footnotesize{$2n+2$} &  \\[0.5 ex]
    \hline\\[-1 ex]

    \footnotesize{\emph{10}} & && &&& & \footnotesize{$(3,6n+1,6n+2,9n+3)$} & && &&& & \footnotesize{$18n+6$} & & &&&& & \footnotesize{$3n+3$} &  \\[0.5 ex]
    \hline\\[-1 ex]

    \footnotesize{\emph{11}} & & &&&& & \footnotesize{$(7,28n-18,42n-27,63n-44)$} & &&& && &  \footnotesize{$126n-81$} & & &&&& & \footnotesize{$7n-1$} &  \\[0.5 ex]
    \hline\\[-1 ex]

    \footnotesize{\emph{12}} & & &&&& & \footnotesize{$(7,28n-17,42n-29,63n-40)$} & &&& & &&\footnotesize{$126n-80$} & && &&& & \footnotesize{$7n+1$} &  \\[0.5 ex]
    \hline\\[-1 ex]

    \footnotesize{\emph{13}} & && &&& & \footnotesize{$(7,28n-13,42n-23,63n-31)$} & &&& && &\footnotesize{$126n-62$} & && &&& & \footnotesize{$7n+2$} &  \\[0.5 ex]
    \hline\\[-1 ex]

    \footnotesize{\emph{14}} & & &&&& & \footnotesize{$(7,28n-10,42n-15,63n-26)$} & &&& && &\footnotesize{$126n-45$} & & &&&&&\footnotesize{$7n+1$} &  \\[0.5 ex]
    \hline\\[-1 ex]

    \footnotesize{\emph{15}} & & &&&& & \footnotesize{$(7,28n-9,42n-17,63n-22)$} && & &&& & \footnotesize{$126n-44$} & & && &&& \footnotesize{$7n+3$} &  \\[0.5 ex]
    \hline\\[-1 ex]

    \footnotesize{\emph{16}} & & &&& && \footnotesize{$(7,28n-6,42n-9,63n-17)$} & &&&& & &\footnotesize{$126n-27$} & & && &&& \footnotesize{$7n+2$} &  \\[0.5 ex]
    \hline\\[-1 ex]

    \footnotesize{\emph{17}} & & &&& &&\footnotesize{$(7,28n-5,42n-11,63n-13)$} & && && &&\footnotesize{$126n-26$} & & &&&& &\footnotesize{$7n+4$} &  \\[0.5 ex]
    \hline\\[-1 ex]

    \footnotesize{\emph{18}} & &&& && &\footnotesize{$(7,28n-2,42n-3,63n-8)$} & & &&&& &\footnotesize{$126n-9$} & && && &&\footnotesize{$7n+3$} &  \\[0.5 ex]
    \hline\\[-1 ex]

    \footnotesize{\emph{19}} & & &&& && \footnotesize{$(7,28n-1,42n-5,63n-4)$} & && &&& &\footnotesize{$126n-8$} & & && &&&\footnotesize{$7n+5$} &  \\[0.5 ex]
    \hline\\[-1 ex]

    \footnotesize{\emph{20}} & &&&&& &\footnotesize{$(7,28n+2,42n+3,63n+1)$} & &&&&& &\footnotesize{$126n+9$} & & &&&& &\footnotesize{$7n+4$} &  \\[0.5 ex]
    \hline\\[-1 ex]

    \footnotesize{\emph{21}} & & &&& && \footnotesize{$(7,28n+3,42n+1,63n+5)$} & & && &&&\footnotesize{$126n+10$} & & &&&& &\footnotesize{$7n+6$} &  \\[0.5 ex]
    \hline\\[-1 ex]

    \footnotesize{\emph{22}} & & &&&& &\footnotesize{$(7,28n+6,42n+9,63n+10)$} & && &&& &\footnotesize{$126n+27$} & & &&& &&\footnotesize{$7n+5$} &  \\[0.5 ex]
    \hline\\[-1 ex]

    \footnotesize{\emph{23}} & &&&& & & \footnotesize{$(2,2n+1,2n+1,4n+1)$} & & &&& &&\footnotesize{$8n+4$} & & && &&& \footnotesize{$1$} &  \\[0.5 ex]
    \hline\\[-1 ex]

    \footnotesize{\emph{24}} & &&& & &&\footnotesize{$(3,3n,3n+1,3n+1)$} & & &&& &&\footnotesize{$9n+3$} & && && &&\footnotesize{$2$} &  \\[0.5 ex]
    \hline\\[-1 ex]

    \footnotesize{\emph{25}} & & &&&& & \footnotesize{$(3,3n+1,3n+2,3n+2)$} & &&&&& &\footnotesize{$9n+6$} & & & &&&& \footnotesize{$2$} &  \\[0.5 ex]
    \hline\\[-1 ex]
    
    \footnotesize{\emph{26}} & & &&&& & \footnotesize{$(3,3n+1,3n+2,6n+1)$} & && &&& &\footnotesize{$12n+5$} & && &&& &\footnotesize{$2$} &  \\[0.5 ex]
    \hline\\[-1 ex]

    \footnotesize{\emph{27}} & && && &&\footnotesize{$(3,3n+1,6n+1,9n)$} & & &&&& & \footnotesize{$18n+3$} & & &&&& &\footnotesize{$2$} &  \\[0.5 ex]
    \hline\\[-1 ex]

    \footnotesize{\emph{28}} & && &&& &\footnotesize{$(3,3n+1,6n+1,9n+3)$} && && & &&\footnotesize{$18n+6$} & && && &&\footnotesize{$2$} &  \\[0.5 ex]
    \hline\\[-1 ex]

    \footnotesize{\emph{29}} & & &&&& &\footnotesize{$(4,2n+3,2n+3,4n+4)$} & && && && \footnotesize{$8n+12$} & & &&& &&\footnotesize{$2$} &  \\[0.5 ex]
    \hline\\[-1 ex]

    \footnotesize{\emph{30}} & && &&& &\footnotesize{$(4,2n+3,4n+6,6n+7)$} & & &&&& &\footnotesize{$12n+18$} & & && &&&\footnotesize{$2$} &  \\[0.5 ex]
    \hline\\[-1 ex]

    \footnotesize{\emph{31}} & &&& && &\footnotesize{$(6,6n+3,6n+5,6n+5)$} & &&&& & &\footnotesize{$18n+15$} & && &&& &\footnotesize{$4$} &  \\[0.5 ex]
    \hline\\[-1 ex]

    \footnotesize{\emph{32}} & & &&&& &\footnotesize{$(6,6n+5,12n+8,18n+9)$} & && &&& & \footnotesize{$36n+24$} & &&&&& &\footnotesize{$4$} &  \\[0.5 ex]
    \hline\\[-1 ex]

    \footnotesize{\emph{33}} & & &&&& & \footnotesize{$(6,6n+5,12n+8,18n+15)$} & &&& && &\footnotesize{$36n+30$} & & &&&& &\footnotesize{$4$} &  \\[0.5 ex]
    \hline\\[-1 ex]

    \footnotesize{\emph{34}} & && &&& &\footnotesize{$(8,4n+5,4n+7,4n+9)$} & &&& && &\footnotesize{$12n+23$} & &&& && &\footnotesize{$6$} &  \\[0.5 ex]
    \hline\\[-1 ex]

    \footnotesize{\emph{35}} & && &&& &\footnotesize{$(9,3n+8,3n+11,6n+13)$} & &&&& & &\footnotesize{$12n+35$} & & &&&& &\footnotesize{$6$} &  \\[0.5 ex]
    \hline\hline
    \end{longtable}
    \end{center}

\bigskip

K-stability of infinite series known so far along with the results in \cite{Johnson2001,Araujo2002,Cheltsov2010,Cheltsov2021a,Boyer2003,Cheltsov2013a,Kim2022,Kim2023} is presented in Table 2.

\begin{center}
    \begin{longtable}{ccccccccccccccccccc}
    \caption{K-stability of Infinite Series}\\
    \hline\\[-1.5 ex]
    \footnotesize{No.} & & &&&& \footnotesize{$(a_0,~I,~\dfrac{a_0}{I})$}&  &&&&& \footnotesize{$\mathrm{K\text{-}unstable}$}& &&&& \footnotesize{$\mathrm{Lower~Indices}$}& \\ [1 ex]
    \hline\hline\\[-1.5 ex]
    \footnotesize{\emph{1}}  & & &&&& \normalsize{$(1,~n,~\frac{1}{n})$}  & &&&&& \normalsize{$n > 2$}& & && &\footnotesize{$
        \text{K-stable}, ~\mathrm{if ~n} \leq 2
        $} &  \\[1 ex]
    \hline\\[-1.5 ex]

    \footnotesize{\emph{2}}  & &&&&& \normalsize{$(1,~n,~\frac{1}{n})$}   & &&&&& \normalsize{$n > 2$}  & & &&&\footnotesize{$\text{K-stable}, ~\mathrm{if ~n} \leq 2 $} & \\[1 ex]
    \hline\\[-1 ex]

    \footnotesize{\emph{3}}  &&&&& & \normalsize{$(1,~n,~\frac{1}{n})$}  &&&&& & \normalsize{$n > 2$}  && &&&\footnotesize{$\text{K-stable}, ~\mathrm{if ~n} \leq 2 $} & \\[1 ex]
    \hline\\[-1 ex]

    \footnotesize{\emph{4}}  &&&&& & \normalsize{$(1,~n,~\frac{1}{n})$} &&&&& & \normalsize{$n > 2$}  & &&&& \footnotesize{$\text{K-stable}, ~\mathrm{if ~n} \leq 2$} & \\[1 ex]
    \hline\\[-1 ex]

    \footnotesize{\emph{5}}  &&&&& & \normalsize{$(1,~n,~\frac{1}{n})$}  &&&&& & \normalsize{$n > 2$}  & &&&& \footnotesize{$\text{K-stable}, ~\mathrm{if ~n} \leq 2$} & \\[1 ex]
    \hline\\[-1 ex]

    \footnotesize{\emph{6}}  &&&&& & \normalsize{$(1,~n,~\frac{1}{n})$}  & &&&&& \normalsize{$n > 1$}  && && &\footnotesize{$
        \text{K-stable}, ~\mathrm{if ~n=1}
        $} &  \\[1 ex]
    \hline\\[-1 ex]

    \footnotesize{\emph{7}}  & &&&&& \normalsize{$(1,~2n,~\frac{1}{2n})$}  &&&&& & \normalsize{$n \geq 1$}  &&&& & \footnotesize{} &  \\[0.9 ex]
    \hline\\[-1 ex]

    \footnotesize{\emph{8}}  &&&&& & \normalsize{$(2,~2n,~\frac{1}{n})$}  & &&&&& \normalsize{$n>3$}  & &&&& \footnotesize{$\left\{
        \begin{array}{ll}
        \text{K-stable},&\mathrm{if ~n=1}\\
        \text{Unknown},&\mathrm{if ~n=2, 3}
        \end{array} \right.$} &  \\[2 ex]
    \hline\\[-1 ex]

    \footnotesize{\emph{9}} &&&&& &  \normalsize{$(2,~2n+2,~\frac{1}{n+1})$}  &&&&& & \normalsize{$n>2$} &&&& &  \footnotesize{$\text{Unknown},\mathrm{~~if ~n} \leq 2$} &  
    \\[1 ex]
    \hline\\[-1 ex]

    \footnotesize{\emph{10}}  & &&&&& \normalsize{$(3,~3n+3,~\frac{1}{n+1})$} &&&& & & \normalsize{$n>2$} &&&&&   \footnotesize{$\text{Unknown},\mathrm{~~if ~n} \leq 2$} &  \\[1 ex]
    \hline\\[-1 ex]

    \footnotesize{\emph{11}}  &&&&& & \normalsize{$(7,~7n-1,~\frac{7}{7n-1})$} &&&& & & \normalsize{$n>3$} & && && \footnotesize{$\text{Unknown}, \mathrm{if ~n} \leq 3$} &  \\[1 ex]
    \hline\\[-1 ex]

    \footnotesize{\emph{12}} & &&&&&  \normalsize{$(7,~7n+1,~\frac{7}{7n+1})$}  &&&&& & \normalsize{$n>2$}  && &&& \footnotesize{$\left\{
        \begin{array}{ll}
        \text{K-stable},&\mathrm{if ~n=1}\\
        \text{Unknown},&\mathrm{if ~n=2}
        \end{array} \right.$} &  \\[2 ex]
    \hline\\[-1 ex]

    \footnotesize{\emph{13}}  &&&&& & \normalsize{$(7,~7n+2,~\frac{7}{7n+2})$}  & &&&&& \normalsize{$n>2$}  &&&& & \footnotesize{$\text{Unknown}, \mathrm{if ~n} \leq 2$} &  \\[1 ex]
    \hline\\[-1 ex]

    \footnotesize{\emph{14}}  &&&&& & \normalsize{$(7,~7n+1,~\frac{7}{7n+1})$} &&&&& &  \normalsize{$n>2$}  && &&& \footnotesize{$\text{Unknown}, \mathrm{if ~n} \leq 2$} &  \\[1 ex]
    \hline\\[-1 ex]

    \footnotesize{\emph{15}}  & &&&&& \normalsize{$(7,~7n+3,~\frac{7}{7n+3})$}  & &&&&& \normalsize{$n>2$}  && &&& \footnotesize{$\text{Unknown}, \mathrm{if ~n} \leq 2$} &  \\[1 ex]
    \hline\\[-1 ex]

    \footnotesize{\emph{16}}  &&&&& & \normalsize{$(7,~7n+2,~\frac{7}{7n+2})$} & & &&&&\normalsize{$n>2$}  && &&& \footnotesize{$\text{Unknown}, \mathrm{if ~n} \leq 2$} &  \\[1 ex]
    \hline\\[-1 ex]

    \footnotesize{\emph{17}} &&&& & &\normalsize{$(7,~7n+4,~\frac{7}{7n+4})$}  &&&&& & \normalsize{$n>2$}  && &&& \footnotesize{$\text{Unknown}, \mathrm{if ~n} \leq 2$} &  \\[1 ex]
    \hline\\[-1 ex]

    \footnotesize{\emph{18}} &&&& & & \normalsize{$(7,~7n+3,~\frac{7}{7n+3})$} &&&&& & \normalsize{$n>2$}  &&& && \footnotesize{$\text{Unknown}, \mathrm{if ~n} \leq 2$} &  \\[1 ex]
    \hline\\[-1 ex]

    \footnotesize{\emph{19}}  &&&&& & \normalsize{$(7,~7n+5,~\frac{7}{7n+5})$}  &&&&& & \normalsize{$n>2$}  & &&& &\footnotesize{$\text{Unknown}, \mathrm{if ~n} \leq 2$} &  \\[1 ex]
    \hline\\[-1 ex]

    \footnotesize{\emph{20}}  &&&&& & \normalsize{$(7,~7n+4,~\frac{7}{7n+4})$} &&&&& & \normalsize{$n>2$}  &&&& & \footnotesize{$\text{Unknown}, \mathrm{if ~n} \leq 2$} &  \\[1 ex]
    \hline\\[-1 ex]

    \footnotesize{\emph{21}}  &&&&& & \normalsize{$(7,~7n+6,~\frac{7}{7n+6})$}  &&&&& & \normalsize{$n>2$}  &&&&&  \footnotesize{$\text{Unknown}, \mathrm{if ~n} \leq 2$} &  \\[1 ex]
    \hline\\[-1 ex]

    \footnotesize{\emph{22}}  & &&&&& \normalsize{$(7,~7n+5,~\frac{7}{7n+5})$}  & & &&&&\normalsize{$n>2$}  && &&& \footnotesize{$\text{Unknown}, \mathrm{if ~n} \leq 2$} &  \\[1 ex]
    \hline\\[-1 ex]

    \footnotesize{\emph{23}}  & &&&&& \small{$(2,~1,~2)$}  &&&&& & \footnotesize{} & &&& &\footnotesize{$\text{K-stable}$} &  \\[1 ex]
    \hline\\[-1 ex]

    \footnotesize{\emph{24}}  &&&&& & \normalsize{$(3,~2,~\frac{3}{2})$}  & &&&&& \footnotesize{} & & && & \footnotesize{$\text{K-stable}$} &  \\[1 ex]
    \hline\\[-1 ex]

    \footnotesize{\emph{25}}  &&&&& &\normalsize{$(3,~2,~\frac{3}{2})$}  &&&&& & \footnotesize{} & &&& & \footnotesize{$\text{K-stable}$} &  \\[1 ex]
    \hline\\[-1 ex]
    
    \footnotesize{\emph{26}}  &&&&& & \normalsize{$(3,~2,~\frac{3}{2})$} &&&&& & \footnotesize{} & & &&& \footnotesize{$\text{K-stable}$} &  \\[1 ex]
    \hline\\[-1 ex]

    \footnotesize{\emph{27}}  & &&&&& \normalsize{$(3,~2,~\frac{3}{2})$} &&&&& & \footnotesize{} & & &&& \footnotesize{$\text{K-stable}$} &  \\[1 ex]
    \hline\\[-1 ex]

    \footnotesize{\emph{28}}  &&&&& & \normalsize{$(3,~2,~\frac{3}{2})$}  &&&&& & \footnotesize{} & & &&& \footnotesize{$\text{K-stable}$} &  \\[1 ex]
    \hline\\[-1 ex]

    \footnotesize{\emph{29}}  &&&&& & \small{$(4,~2,~2)$} &  &&&&& \footnotesize{}  &&&& &\footnotesize{$\text{K-stable}$} &  \\[1 ex]
    \hline\\[-1 ex]

    \footnotesize{\emph{30}}  & &&&&& \small{$(4,~2,~2)$}  &&&&& & \footnotesize{}  &&&& &\footnotesize{$\text{K-stable}$} &  \\[1 ex]
    \hline\\[-1 ex]

    \footnotesize{\emph{31}}  &&&&& & \normalsize{$(6,~4,~\frac{3}{2})$} &&&&& & \footnotesize{} && &&&  \footnotesize{$\text{K-stable}$} &  \\[1 ex]
    \hline\\[-1 ex]

    \footnotesize{\emph{32}}  &&&&& & \normalsize{$(6,~4,~\frac{3}{2})$}  &&&&& & \footnotesize{} & &&&&  \footnotesize{$\text{K-stable}$} &  \\[1 ex]
    \hline\\[-1 ex]

    \footnotesize{\emph{33}}  &&&&& & \normalsize{$(6,~4,~\frac{3}{2})$}  &&&&& & \footnotesize{} & & &&& \footnotesize{$\text{K-stable}$} &  \\[1 ex]
    \hline\\[-1 ex]

    \footnotesize{\emph{34}} & &&&& & \normalsize{$(8,~6,~\frac{4}{3})$} &&& && & \footnotesize{} & & &&& \footnotesize{$\text{K-stable}$} &  \\[1 ex]
    \hline\\[-1 ex]

    \footnotesize{\emph{35}} & &&&& & \normalsize{$(9,~6,~\frac{3}{2})$} & &&&& & \footnotesize{} & & &&& \footnotesize{$\text{K-stable}$} &  \\[1 ex]
    \hline\hline
    \end{longtable}
    \end{center}

    \bigskip

    Each linear system and singular locus is presented in Table 3 to prove that a non-lc point is contained in the section $H_x$ defined by $x=0$ as in Lemma \ref*{Lemma:lc for S_k-H_x}. In the table, $\mathrm{deg}(\mathscr{M})$ means the degree of a generic member of the linear system $\mathscr{M}$.

\begin{center}
    \begin{longtable}{cccccccccc}
        \caption{Singular Points and Linear Systems for Non-lc Locus}\\
    \hline\\[-1.5 ex]
    \footnotesize{No.}  & & & \footnotesize{Singular Locus}& & &  \footnotesize{$\mathscr{M}$}& && \footnotesize{$\mathrm{deg}(\mathscr{M})$} \\ [1 ex]
    \hline\hline\\[-1.5 ex]
    \footnotesize{\emph{1}}  & & & \footnotesize{$\mathsf{p}_{_y}, ~\mathsf{p}_{_t}$} &  && \footnotesize{$\Big \vert ~\alpha x^{4n-3} + \beta z + \gamma x^{n-1}y =0 ~\Big \vert$}& & & \footnotesize{$4n-3$}   \\[1.5 ex]
    \hline\\[-1.5 ex]

    \footnotesize{\emph{2}} & & & \footnotesize{$\mathsf{p}_{_z}, ~\mathsf{p}_{_{3n-2}}=[0:a:0:b] \cap S_2$} & &&  \footnotesize{$\Big \vert ~  ~\alpha x^{4n-3} + \beta z +\gamma x^{n-1}y =0 ~\Big \vert$}  & & & \footnotesize{$4n-3$}  \\[1.5 ex]
    \hline\\[-1 ex]

    \footnotesize{\emph{3}} & & & \footnotesize{$\mathsf{p}_{_y}, ~\mathsf{p}_{_z}$} & &&  \footnotesize{$\Big \vert ~ \alpha x^{6n-5}+\beta z+\gamma x^{2n-2}y =0 ~\Big \vert$} & & & \footnotesize{$6n-5$}  \\[1.5 ex]
    \hline\\[-1 ex]

    \footnotesize{\emph{4}} & & & \footnotesize{$\mathsf{p}_{_y}, ~\mathsf{p}_{_{5n-4}}=[0:0:a:b] \cap S_4$} &&  & \footnotesize{$\Big \vert ~\alpha x^{10n-8}+\beta z +\gamma x^{4n-3}y=0 ~\Big \vert$} & & & \footnotesize{$10n-8$}  \\[1.5 ex]
    \hline\\[-1 ex]

    \footnotesize{\emph{5}} & & & \footnotesize{$\mathsf{p}_{_z}, ~\mathsf{p}_{_{3n-2}}=[0:a:0:b]\cap S_5$} & &&  \footnotesize{$\Big \vert ~\alpha x^{10n-7}+\beta z+\gamma x^{4n-3}y=0 ~\Big \vert$} & & & \footnotesize{$10n-7$}  \\[1.5 ex]
    \hline\\[-1 ex]

    \footnotesize{\emph{6}} & & & \footnotesize{$\mathsf{p}_{_t}, ~\mathsf{p}_{_{2n-1}}=[0:a:b:0]\cap S_6$} & & & \footnotesize{$\Big \vert ~\alpha x^{10n-5}+\beta z+\gamma x^{4n-2}y=0 ~\Big \vert$} & & & \footnotesize{$10n-5$}   \\[1.5 ex]
    \hline\\[-1 ex]

    \footnotesize{\emph{7}} & & & \footnotesize{$\mathsf{p}_{_y}, ~\mathsf{p}_{_t}, ~\mathsf{p}_{_{4n-1}}=[0:a:b:0]\cap S_7$} & & & \footnotesize{$\Big \vert ~\alpha x^{12n-3}+\beta z +\gamma x^{4n-1}y=0~\Big \vert$} & & & \footnotesize{$12n-3$}   \\[1.5 ex]
    \hline\\[-1 ex]

    \footnotesize{\emph{8}} & & & \footnotesize{$
        \begin{array}{rl}
            &\mathsf{p}_{_t}, ~\mathsf{p}_{_2}=[a:0:b:0]\cap S_8, \\ 
            & \hspace{1 cm}~\mathsf{p}_{_{2n-1}}=[0:a:b:0]\cap S_8
        \end{array} $} & &&  \footnotesize{$\Big \vert ~\alpha x^{6n-3}+\beta y^2+\gamma x^{2n-1}z=0 ~\Big \vert$} & & & \footnotesize{$12n-6$}   \\[1.5 ex]
    \hline\\[-1 ex]

    \footnotesize{\emph{9}} & & & \footnotesize{$
    \begin{array}{rl}
    &\mathsf{p}_{_y}, ~\mathsf{p}_{_2}=[a:0:b:0]\cap S_9,\\ 
    & \hspace{1 cm}~\mathsf{p}_{_{4n+1}}=[0:0:a:b]\cap S_9
    \end{array} $} & & & \footnotesize{$\Big \vert ~\alpha x^{6n+1}+\beta y^2+\gamma x^{2n}z=0 ~ \Big \vert$} & & & \footnotesize{$12n+2$}   \\[1.5 ex]
    \hline\\[-1 ex]

    \footnotesize{\emph{10}} & & & \footnotesize{$
    \begin{array}{rl}
    &\mathsf{p}_{_y}, ~\mathsf{p}_{_3}=[a:0:0:b]\cap S_{10}, \\ 
    & \hspace{1 cm}~\mathsf{p}_{_{3n+1}}=[0:0:a:b]\cap S_{10}
    \end{array} $} & & & \footnotesize{$\Big \vert ~\alpha x^{3n+1}+\beta t =0 ~\Big \vert $} & & & \footnotesize{$9n+3$}   \\[1.5 ex]
    \hline\\[-1 ex]

    \footnotesize{\emph{11}} & & & \footnotesize{$\mathsf{p}_{_x}, ~\mathsf{p}_{_y}, ~\mathsf{p}_{_t}, ~\mathsf{p}_{_{14n-9}}=[0:a:b:0]\cap S_{11}$} & & & \footnotesize{$\Big \vert ~\alpha x^{28n-18}+\beta y^7 =0 ~\Big \vert $} & & & \footnotesize{$196n-126$}   \\[2ex]
    \hline\\[-1 ex]

    \footnotesize{\emph{12}} & & & \footnotesize{$\mathsf{p}_{_x}, ~\mathsf{p}_{_y}, ~\mathsf{p}_{_z}$} & & & \footnotesize{$\Big \vert ~\alpha x^{28n-17}+\beta y^7 =0 ~\Big \vert $} & & & \footnotesize{$196n-119$}   \\[2 ex]
    \hline\\[-1 ex]

    \footnotesize{\emph{13}}  & & & \footnotesize{$\mathsf{p}_{_x}, ~\mathsf{p}_{_y}, ~\mathsf{p}_{_z}$} & &&  \footnotesize{$\Big \vert ~\alpha x^{28n-13}+\beta y^7 =0 ~\Big \vert$} & & & \footnotesize{$196n-91$}   \\[2 ex]
    \hline\\[-1 ex]

    \footnotesize{\emph{14}} & & & \footnotesize{$\mathsf{p}_{_x}, ~\mathsf{p}_{_y}, ~\mathsf{p}_{_t}, ~\mathsf{p}_{_{14n-5}}=[0:a:b:0] \cap S_{14}$} & & & \footnotesize{$\Big \vert ~\alpha x^{28n-10}+\beta y^7 =0 ~\Big \vert$} & & & \footnotesize{$196n-70$}   \\[2 ex]
    \hline\\[-1 ex]

    \footnotesize{\emph{15}} & & & \footnotesize{$\mathsf{p}_{_x}, ~\mathsf{p}_{_y}, ~\mathsf{p}_{_z}$} & & & \footnotesize{$\Big \vert ~\alpha x^{28n-9}+\beta y^7 =0 ~\Big \vert$} & & & \footnotesize{$196n-63$}   \\[2 ex]
    \hline\\[-1 ex]

    \footnotesize{\emph{16}} & & & \footnotesize{$\mathsf{p}_{_x}, ~\mathsf{p}_{_y}, ~\mathsf{p}_{_t}, ~\mathsf{p}_{_{14n-3}}=[0:a:b:0] \cap S_{16}$} & &&  \footnotesize{$\Big \vert ~\alpha x^{28n-6}+\beta y^7 =0 ~\Big \vert$} & & & \footnotesize{$196n-42$}   \\[2 ex]
    \hline\\[-1 ex]

    \footnotesize{\emph{17}} & & & \footnotesize{$\mathsf{p}_{_x}, ~\mathsf{p}_{_y}, ~\mathsf{p}_{_z}$} & &&  \footnotesize{$\Big \vert ~\alpha x^{28n-5}+\beta y^7 =0 ~\Big \vert$} & & & \footnotesize{$196n-35$}   \\[2 ex]
    \hline\\[-1 ex]

    \footnotesize{\emph{18}} & & & \footnotesize{$\mathsf{p}_{_x}, ~\mathsf{p}_{_y}, ~\mathsf{p}_{_t}, ~\mathsf{p}_{_{14n-1}}=[0:a:b:0] \cap S_{18}$} & &&  \footnotesize{$\Big \vert ~\alpha x^{28n-2}+\beta y^7 =0 ~\Big \vert$} & & & \footnotesize{$196n-14$}   \\[2 ex]
    \hline\\[-1 ex]

    \footnotesize{\emph{19}} & & & \footnotesize{$\mathsf{p}_{_x}, ~\mathsf{p}_{_y}, ~\mathsf{p}_{_z}$} & &&  \footnotesize{$\Big \vert ~\alpha x^{28n-1}+\beta y^7 =0 ~\Big \vert$} & & & \footnotesize{$196n-7$}   \\[2 ex]
    \hline\\[-1 ex]

    \footnotesize{\emph{20}} & & & \footnotesize{$\mathsf{p}_{_x}, ~\mathsf{p}_{_y}, ~\mathsf{p}_{_t}, ~\mathsf{p}_{_{14n+1}}=[0:a:b:0] \cap S_{20}$} && &  \footnotesize{$\Big \vert ~\alpha x^{28n+2}+\beta y^7 =0 ~\Big \vert$} & & & \footnotesize{$196n+14$}   \\[2 ex]
    \hline\\[-1 ex]

    \footnotesize{\emph{21}} & & & \footnotesize{$\mathsf{p}_{_x}, ~\mathsf{p}_{_y}, ~\mathsf{p}_{_z}$} & &&  \footnotesize{$\Big \vert ~\alpha x^{28n+3}+\beta y^7 =0 ~\Big \vert$} & & & \footnotesize{$196n+21$}   \\[2 ex]
    \hline\\[-1 ex]

    \footnotesize{\emph{22}} & & & \footnotesize{$\mathsf{p}_{_x}, ~\mathsf{p}_{_y}, ~\mathsf{p}_{_t}, ~\mathsf{p}_{_{14n+3}}=[0:a:b:0] \cap S_{22}$} & &&  \footnotesize{$\Big \vert ~\alpha x^{28n+6}+\beta y^7 =0 ~\Big \vert$} & & & \footnotesize{$196n+42$}   \\[1.5 ex]    
    \hline\hline
    \end{longtable}
    \end{center}

    \bigskip

Each $D^{\ast}$ to prove that the base locus of the corresponding pencil cannot be empty as in Lemma \ref*{lem:baseptfree} is presented in Table 4. 

\begin{center}
    \begin{longtable}{ccccccccccccccccc}
        \caption{Non-emptiness of Base Locus of the Pencil induced by Cylinder}\\
    \hline\\[-1.5 ex]
    \footnotesize{No.} $\hspace{2.5 cm}$ & & & &&&\footnotesize{$D^{\ast}$} $\hspace{2 cm}$& & &&&&  \\ [1 ex]
    \hline\hline\\[-1.5 ex]
    \footnotesize{\emph{1 - 6, 8}} $\hspace{2.5 cm}$& &&&& & \footnotesize{$\dfrac{n}{n-2}D-\dfrac{2n}{n-2}H_x$} $\hspace{3 cm}$& &&&& &  \\[2ex]
    \hline\\[-1.5 ex]

    \footnotesize{\emph{7}}$\hspace{2.5 cm}$ &&&& & & \footnotesize{$\dfrac{n}{n-1}D-\dfrac{2n}{n-1}H_x$} $\hspace{3 cm}$& & &&&&   \\[2 ex]
    \hline\\[-1 ex]

    \footnotesize{\emph{9, 10}} $\hspace{2.5 cm}$&&&& & & \footnotesize{$\dfrac{n+1}{n-1}D-\dfrac{2(n+1)}{n-1}H_x$} $\hspace{3 cm}$& &&&& &   \\[2 ex]
    \hline\\[-1 ex]

    \footnotesize{\emph{11}} $\hspace{2.5 cm}$& &&&& & \footnotesize{$
        \dfrac{(7n-1)}{(7n-15)}D-\dfrac{2(7n-1)}{(7n-15)}H_x$} $\hspace{3 cm}$& &&&& &  \\[2ex]
    \hline\\[-1 ex]

    \footnotesize{\emph{12, 14}} $\hspace{2.5 cm}$& &&&& & \footnotesize{$\dfrac{(7n+1)}{(7n-13)}D-\dfrac{2(7n+1)}{(7n-13)}H_x$}
    $\hspace{3 cm}$& &&&& &  \\[2 ex]
    \hline\\[-1 ex]

    \footnotesize{\emph{13, 16}} $\hspace{2.5 cm}$& & &&&& \footnotesize{$\dfrac{(7n+2)}{(7n-12)}D-\dfrac{2(7n+2)}{(7n-12)}H_x
       $} $\hspace{3 cm}$& &&&& &  \\[2 ex]
    \hline\\[-1 ex]

    \footnotesize{\emph{15, 18}} $\hspace{2.5 cm}$& & &&&& \footnotesize{$\dfrac{(7n+3)}{(7n-11)}D-\dfrac{2(7n+3)}{(7n-11)}H_x
        $} $\hspace{3 cm}$& &&&& &  \\[2 ex]
    \hline\\[-1 ex]

    \footnotesize{\emph{17, 20}} $\hspace{2.5 cm}$& & &&&& \footnotesize{$\dfrac{(7n+4)}{(7n-10)}D-\dfrac{2(7n+4)}{(7n-10)}H_x
        $} $\hspace{3 cm}$& & &&&&  \\[2 ex]
    \hline\\[-1 ex]

    \footnotesize{\emph{19, 22}} $\hspace{2.5 cm}$& & &&&& \footnotesize{$\dfrac{(7n+5)}{(7n-9)}D-\dfrac{2(7n+5)}{(7n-9)}H_x
        $} $\hspace{3 cm}$& & &&&&   \\[2 ex]
    \hline\\[-1 ex]

    \footnotesize{\emph{21}} $\hspace{2.5 cm}$& & &&&& \footnotesize{$\dfrac{(7n+6)}{(7n-8)}D-\dfrac{2(7n+6)}{(7n-8)}H_x
        $} $\hspace{3 cm}$& & &&&&   \\[2ex]
    \hline\hline
    \end{longtable}
    \end{center}

    \bigskip

Each absurd inequality where $H_x$ is irreducible as in Theorem \ref*{main1} is presented in Table 5.

\begin{center}
    \begin{longtable}{cccccc}
        \caption{Absurd Inequalities for Irreducible $H_x$}\\
    \hline\\[-1.5 ex]
    \footnotesize{No.} & & & \footnotesize{Absurd Inequalities} &  &\hspace{-0.3 cm} \\ [1 ex]
    \hline\hline\\[-1.5 ex]
    \footnotesize{\emph{1}} && &  \footnotesize{$\dfrac{3n}{(3n-2)(6n-5)}=H_x \cdot D \geq \mult_{\mathsf{p}}(H_x)\cdot \mult_{\mathsf{p}}(D)> \left\{
        \begin{array}{ll}
        1, ~~\mathrm{if} ~\mathsf{p} ~~\mathrm{is ~~smooth}\\\\   
        \dfrac{2}{3n-2}, ~~\mathrm{if} ~\mathsf{p}=\mathsf{p}_{_y}\\\\
        \dfrac{2}{6n-5}, ~~\mathrm{if} ~\mathsf{p}=\mathsf{p}_{_t}\\
        \end{array} \right. $} &  & \hspace{-0.3 cm}  \\[8ex]
    \hline\\[-1.5 ex]

    \footnotesize{\emph{3}} & &&  \footnotesize{$\dfrac{2n}{(4n-3)(6n-5)}=H_x \cdot D \geq \mult_{\mathsf{p}}(H_x)\cdot \mult_{\mathsf{p}}(D)> \left\{
        \begin{array}{ll}
        1, ~~\mathrm{if} ~\mathsf{p} ~~\mathrm{is ~~smooth}\\\\   
        \dfrac{2}{4n-3}, ~~\mathrm{if} ~\mathsf{p}=\mathsf{p}_{_y}\\\\
        \dfrac{2}{6n-5}, ~~\mathrm{if} ~\mathsf{p}=\mathsf{p}_{_z}\\
        \end{array} \right. $} & & \hspace{-0.3 cm} \\[8 ex]
    \hline\\[-1 ex]

    \footnotesize{\emph{4}} & & & \footnotesize{$\dfrac{n}{(5n-4)(6n-5)}=T\cdot D \geq \mult_{\mathsf{p}}(T)\cdot \mult_{\mathsf{p}}(D)> \left\{
        \begin{array}{ll}
        1, ~~\mathrm{if} ~\mathsf{p} ~~\mathrm{is ~~smooth}\\\\   
        \dfrac{2}{6n-5}, ~~\mathrm{if} ~\mathsf{p}=\mathsf{p}_{_y}\\\\
        \dfrac{2}{5n-4}, ~~\mathrm{if} ~\mathsf{p}=\mathsf{p}_{_{5n-4}}\\
        \end{array} \right. $} & & \hspace{-0.3 cm}   \\[8ex]
    \hline\\[-1 ex]

    \footnotesize{\emph{5}} & & & \footnotesize{$\dfrac{n}{(3n-2)(10n-7)}=H_x \cdot D \geq \mult_{\mathsf{p}}(H_x)\cdot \mult_{\mathsf{p}}(D)> \left\{
        \begin{array}{ll}
        1, ~~\mathrm{if} ~\mathsf{p} ~~\mathrm{is ~~smooth}\\\\   
        \dfrac{2}{10n-7}, ~~\mathrm{if} ~\mathsf{p}=\mathsf{p}_{_z}\\\\
        \dfrac{2}{3n-2}, ~~\mathrm{if} ~\mathsf{p}=\mathsf{p}_{_{3n-2}}\\
        \end{array} \right.$} & & \hspace{-0.3 cm} \\[8 ex]
    \hline\\[-1 ex]

    \footnotesize{\emph{6}} & & & \footnotesize{$\dfrac{n}{(2n-1)(15n-8)}=H_x \cdot D \geq \mult_{\mathsf{p}}(H_x)\cdot \mult_{\mathsf{p}}(D)> \left\{
        \begin{array}{ll}
        1, ~~\mathrm{if} ~\mathsf{p} ~~\mathrm{is ~~smooth}\\\\   
        \dfrac{2}{15n-8}, ~~\mathrm{if} ~\mathsf{p}=\mathsf{p}_{_t}\\\\
        \dfrac{2}{2n-1}, ~~\mathrm{if} ~\mathsf{p}=\mathsf{p}_{_{2n-1}}\\
        \end{array} \right.$} & & \hspace{-0.3 cm}  \\[8 ex]
    \hline\\[-1 ex]

    \footnotesize{\emph{8}} &  && \footnotesize{$\dfrac{2n}{(2n-1)(12n-7)}=H_x \cdot D \geq \mult_{\mathsf{p}}(H_x)\cdot \mult_{\mathsf{p}}(D)> \left\{
        \begin{array}{ll}
        1, ~~\mathrm{if ~~either} ~\mathsf{p} ~~\mathrm{is ~~smooth ~~or} ~~\mathsf{p}=\mathsf{p}_{_2}\\\\   
        \dfrac{2}{12n-7}, ~~\mathrm{if} ~\mathsf{p}=\mathsf{p}_{_t}\\\\
        \dfrac{2}{2n-1}, ~~\mathrm{if} ~\mathsf{p}=\mathsf{p}_{_{2n-1}}\\
        \end{array} \right.$} &  & \hspace{-0.3 cm}  \\[8 ex]
    \hline\\[-1 ex]

    \footnotesize{\emph{9}} & &  & \footnotesize{$\dfrac{2n+2}{(4n+1)(6n+1)}=H_x \cdot D \geq \mult_{\mathsf{p}}(H_x)\cdot \mult_{\mathsf{p}}(D)> \left\{
        \begin{array}{ll}
        1, ~~\mathrm{if ~~either} ~\mathsf{p} ~~\mathrm{is ~~smooth ~~or} ~~\mathsf{p}=\mathsf{p}_{_2}\\\\ 
        \dfrac{2}{6n+1}, ~~\mathrm{if} ~\mathsf{p}=\mathsf{p}_{_y}\\\\
        \dfrac{2}{4n+1}, ~~\mathrm{if} ~\mathsf{p}=\mathsf{p}_{_{4n+1}}\\
        \end{array} \right.$} &  & \hspace{-0.3 cm}  \\[8 ex]
    \hline\\[-1 ex]

    \footnotesize{\emph{10}} & & & \footnotesize{$\dfrac{3n+3}{(3n+1)(6n+1)}=H_x\cdot D \geq \mult_{\mathsf{p}}(H_x)\cdot \mult_{\mathsf{p}}(D)> \left\{
        \begin{array}{ll}
        1, ~~\mathrm{if} ~\mathsf{p} ~~\mathrm{is ~~smooth}\\\\   
        \dfrac{2}{6n+1}, ~~\mathrm{if} ~\mathsf{p}=\mathsf{p}_{_y}\\\\
        \dfrac{2}{3}, ~~\mathrm{if} ~\mathsf{p}=\mathsf{p}_{_3}\\\\
        \dfrac{2}{3n+1}, ~~\mathrm{if} ~\mathsf{p}=\mathsf{p}_{_{3n+1}}\\
        \end{array} \right.$} & & \hspace{-0.3 cm}   \\[11 ex]
    \hline\\[-1 ex]

    \footnotesize{\emph{12}} & & & \footnotesize{$
        \dfrac{2(7n+1)}{(28n-17)(42n-29)}=H_x\cdot D \geq \mult_{\mathsf{p}}(H_x)\cdot \mult_{\mathsf{p}}(D)>\left\{
            \begin{array}{ll}
            1, ~~\mathrm{if} ~\mathsf{p} ~~\mathrm{is ~~smooth}\\\\   
            \dfrac{2}{28n-17}, ~~\mathrm{if} ~\mathsf{p}=\mathsf{p}_{_y}\\\\
            \dfrac{2}{42n-29}, ~~\mathrm{if} ~\mathsf{p}=\mathsf{p}_{_z}\\
            \end{array} \right.
        $} &  & \hspace{-0.3 cm} \\[8 ex]
    \hline\\[-1 ex]

    \footnotesize{\emph{13}} & & & \footnotesize{$\dfrac{2(7n+2)}{(28n-13)(42n-23)}=H_x\cdot D \geq \mult_{\mathsf{p}}(H_x)\cdot \mult_{\mathsf{p}}(D)>\left\{
            \begin{array}{ll}
            1, ~~\mathrm{if} ~\mathsf{p} ~~\mathrm{is ~~smooth}\\\\   
            \dfrac{2}{28n-13}, ~~\mathrm{if} ~\mathsf{p}=\mathsf{p}_{_y}\\\\
            \dfrac{2}{42n-23}, ~~\mathrm{if} ~\mathsf{p}=\mathsf{p}_{_z}\\
            \end{array} \right.$} & & \hspace{-0.3 cm}  \\[8 ex]
    \hline\\[-1 ex]

    \footnotesize{\emph{15}}  && & \footnotesize{$\dfrac{2(7n+3)}{(28n-9)(42n-17)}=H_x\cdot D \geq \mult_{\mathsf{p}}(H_x)\cdot \mult_{\mathsf{p}}(D)>\left\{
            \begin{array}{ll}
            1, ~~\mathrm{if} ~\mathsf{p} ~~\mathrm{is ~~smooth}\\\\   
            \dfrac{2}{28n-9}, ~~\mathrm{if} ~\mathsf{p}=\mathsf{p}_{_y}\\\\
            \dfrac{2}{42n-17}, ~~\mathrm{if} ~\mathsf{p}=\mathsf{p}_{_z}\\
            \end{array} \right.$} &  & \hspace{-0.3 cm}  \\[8 ex]
    \hline\\[-1 ex]

    \footnotesize{\emph{17}} &&  & \footnotesize{$\dfrac{2(7n+4)}{(28n-5)(42n-11)}=H_x\cdot D \geq \mult_{\mathsf{p}}(H_x)\cdot \mult_{\mathsf{p}}(D)>\left\{
            \begin{array}{ll}
            1, ~~\mathrm{if} ~\mathsf{p} ~~\mathrm{is ~~smooth}\\\\   
            \dfrac{2}{28n-5}, ~~\mathrm{if} ~\mathsf{p}=\mathsf{p}_{_y}\\\\
            \dfrac{2}{42n-11}, ~~\mathrm{if} ~\mathsf{p}=\mathsf{p}_{_z}\\
            \end{array} \right.$} &  &\hspace{-0.3 cm}  \\[8 ex]
    \hline\\[-1 ex]

    \footnotesize{\emph{19}} & && \footnotesize{$\dfrac{2(7n+5)}{(28n-1)(42n-5)}=H_x\cdot D \geq \mult_{\mathsf{p}}(H_x)\cdot \mult_{\mathsf{p}}(D)>\left\{
            \begin{array}{ll}
            1, ~~\mathrm{if} ~\mathsf{p} ~~\mathrm{is ~~smooth}\\\\   
            \dfrac{2}{28n-1}, ~~\mathrm{if} ~\mathsf{p}=\mathsf{p}_{_y}\\\\
            \dfrac{2}{42n-5}, ~~\mathrm{if} ~\mathsf{p}=\mathsf{p}_{_z}\\
            \end{array} \right.$} &  & \hspace{-0.3 cm}  \\[8 ex]
    \hline\\[-1 ex]

    \footnotesize{\emph{21}} &  && \footnotesize{$\dfrac{2(7n+6)}{(28n+3)(42n+1)}=H_x\cdot D \geq \mult_{\mathsf{p}}(H_x)\cdot \mult_{\mathsf{p}}(D)>\left\{
            \begin{array}{ll}
            1, ~~\mathrm{if} ~\mathsf{p} ~~\mathrm{is ~~smooth}\\\\   
            \dfrac{2}{28n+3}, ~~\mathrm{if} ~\mathsf{p}=\mathsf{p}_{_y}\\\\
            \dfrac{2}{42n+1}, ~~\mathrm{if} ~\mathsf{p}=\mathsf{p}_{_z}\\
            \end{array} \right.$} &  & \hspace{-0.3 cm}  \\[8ex]
    \hline\hline
    \end{longtable}
    \end{center}

\newpage

The absurd inequalities, \emph{(A0)} through \emph{(A3)}, where $H_x$ is reducible as in Theorem \ref*{main1} are presented in Table 6.

    \begin{center}
        \begin{longtable}{ccc}
            \caption{Absurd Inequalities for Reducible $H_x$}\\
        \hline\\[-1.5 ex]
        \footnotesize{No.} \hspace{-1cm} &   \footnotesize{Absurd Inequalities} &  \hspace{-1cm} \\ [1 ex]
        \hline\hline\\[-1.5 ex]
    
        \footnotesize{\emph{7}}  \hspace{-1cm} & \footnotesize{$\begin{array}{ll}
        &\emph{(A0)} \begin{array}{ll}
                L\cdot(&\hspace{-0.3 cm}-K_{S_7})=\dfrac{n}{(4n-1)(18n-5)}, ~~R\cdot(-K_{S_7})=\dfrac{2n}{(4n-1)(18n-5)}, ~~L\cdot R=\dfrac{3}{18n-5}\\\\
                &L^2=\dfrac{1}{(8n-2)(18n-5)}-\dfrac{3}{18n-5}, ~~R^2=\dfrac{2}{(8n-2)(18n-5)}-\dfrac{3}{18n-5}\\\\
               \end{array} \\\\
        &\emph{(A1)} ~~\dfrac{3n}{(4n-1)(18n-5)}=H_x\cdot D \geq \mult_{\mathsf{p}}(H_x)\cdot \mult_{\mathsf{p}}(D)> \left\{
              \begin{array}{ll}
              1, ~~\mathrm{if} ~\mathsf{p} ~~\mathrm{is ~~smooth}\\\\   
              \dfrac{2}{8n-2}, ~~\mathrm{if} ~\mathsf{p}=\mathsf{p}_{_y}\\\\
              \dfrac{2}{18n-5}, ~~\mathrm{if} ~\mathsf{p}=\mathsf{p}_{_t}\\\\
              \dfrac{2}{4n-1}, ~~\mathrm{if} ~\mathsf{p}=\mathsf{p}_{_{4n-1}}\\\\
              \end{array} \right.\\\\

        &\emph{(A2)} ~~\dfrac{n}{(4n-1)(18n-5)}=L\cdot D > \dfrac{1}{18n-5}
        ~~~\left(\mathrm{resp}. ~~\dfrac{2n}{(4n-1)(18n-5)}=R\cdot D > \dfrac{1}{18n-5} \right)\\\\\\

        &\emph{(A3)} ~~\dfrac{2n+\lambda_{R}(12n-4)}{(4n-1)(18n-5)}=R\cdot (D-\lambda_{R}R)= R \cdot\Delta > \left\{ 
                    \begin{array}{ll}
                    1, ~~\mathrm{if} ~\mathsf{p} ~~\mathrm{is ~~smooth}\\\\   
                    \dfrac{1}{4n-1}, ~~\mathrm{if} ~\mathsf{p}=\mathsf{p}_{_{4n-1}}\\
                    \end{array} \right.    
            \end{array}$} &  \hspace{-1cm} \\[30ex]
        \hline\\[-1 ex]

        \footnotesize{\emph{11}}  \hspace{-1cm} & \footnotesize{$\begin{array}{ll}
            &\emph{(A0)} \begin{array}{ll}
                    L\cdot(&\hspace{-0.3 cm}-K_{S_{11}})=\dfrac{7n-1}{(28n-18)(63n-44)}, ~~R\cdot(-K_{S_{11}})=\dfrac{14n-2}{(28n-18)(63n-44)}, ~~L\cdot R=\dfrac{3}{63n-44}\\\\
                    &L^2=\dfrac{7}{(28n-18)(63n-44)}-\dfrac{3}{63n-44}, ~~R^2=\dfrac{14}{(28n-18)(63n-44)}-\dfrac{3}{63n-44}\\\\
                   \end{array} \\\\
            &\emph{(A1)} ~~\dfrac{3(7n-1)}{(28n-18)(63n-44)}=H_x\cdot D \geq \mult_{\mathsf{p}}(H_x)\cdot \mult_{\mathsf{p}}(D)>\left\{
                \begin{array}{ll}
                1, ~~\mathrm{if} ~\mathsf{p} ~~\mathrm{is ~~smooth}\\\\   
                \dfrac{2}{28n-18}, ~~\mathrm{if} ~\mathsf{p}=\mathsf{p}_{_y}\\\\
                \dfrac{2}{63n-44}, ~~\mathrm{if} ~\mathsf{p}=\mathsf{p}_{_t}\\\\
                \dfrac{2}{14n-9}, ~~\mathrm{if} ~\mathsf{p}=\mathsf{p}_{_{14n-9}}\\
                \end{array} \right.\\\\
            &\emph{(A2)} ~~\dfrac{7n-1}{(28n-18)(63n-44)}=L\cdot D > \dfrac{1}{63n-44}
            ~~~\left(\mathrm{resp}. ~~\dfrac{14n-2}{(28n-18)(63n-44)}=R\cdot D > \dfrac{1}{63n-44} \right)\\\\\\
            
            &\emph{(A3)} ~~\dfrac{14n-2+\lambda_{R}(84n-68)}{(28n-18)(63n-44)}=R\cdot (D-\lambda_{R}R)= R \cdot\Delta > \left\{ 
                        \begin{array}{ll}
                        1, ~~\mathrm{if} ~\mathsf{p} ~~\mathrm{is ~~smooth}\\\\   
                        \dfrac{1}{14n-9}, ~~\mathrm{if} ~\mathsf{p}=\mathsf{p}_{_{14n-9}}\\
                        \end{array} \right.    
                \end{array}$} & \hspace{-1cm}  \\[30ex]
                \hline\\[-1 ex]

        \footnotesize{\emph{14}}  \hspace{-1cm} & \footnotesize{$\begin{array}{ll}
            &\emph{(A0)} \begin{array}{ll}
                    L\cdot(&\hspace{-0.3 cm}-K_{S_{14}})=\dfrac{7n+1}{(28n-10)(63n-26)}, ~~R\cdot(-K_{S_{14}})=\dfrac{14n+2}{(28n-10)(63n-26)}, ~~L\cdot R=\dfrac{3}{63n-26}\\\\
                    &L^2=\dfrac{7}{(28n-10)(63n-26)}-\dfrac{3}{63n-26}, ~~R^2=\dfrac{14}{(28n-10)(63n-26)}-\dfrac{3}{63n-26}\\\\
                   \end{array} \\\\
            &\emph{(A1)} ~~\dfrac{3(7n+1)}{(28n-10)(63n-26)}=H_x \cdot D \geq \mult_{\mathsf{p}}(H_x)\cdot \mult_{\mathsf{p}}(D)>\left\{
                \begin{array}{ll}
                1, ~~\mathrm{if} ~\mathsf{p} ~~\mathrm{is ~~smooth}\\\\   
                \dfrac{2}{28n-10}, ~~\mathrm{if} ~\mathsf{p}=\mathsf{p}_{_y}\\\\
                \dfrac{2}{63n-26}, ~~\mathrm{if} ~\mathsf{p}=\mathsf{p}_{_t}\\\\
                \dfrac{2}{14n-5}, ~~\mathrm{if} ~\mathsf{p}=\mathsf{p}_{_{14n-5}}\\
                \end{array} \right.\\\\
            &\emph{(A2)} ~~\dfrac{7n+1}{(28n-10)(63n-26)}=L\cdot D > \dfrac{1}{63n-26}
            ~~~\left(\mathrm{resp}. ~~\dfrac{14n+2}{(28n-10)(63n-26)}=R\cdot D > \dfrac{1}{63n-26} \right)\\\\\\
            
            &\emph{(A3)} ~~\dfrac{14n+2+\lambda_{R}(84n-44)}{(28n-10)(63n-26)}=R\cdot (D-\lambda_{R}R)= R \cdot\Delta > \left\{ 
                        \begin{array}{ll}
                        1, ~~\mathrm{if} ~\mathsf{p} ~~\mathrm{is ~~smooth}\\\\   
                        \dfrac{1}{14n-5}, ~~\mathrm{if} ~\mathsf{p}=\mathsf{p}_{_{14n-5}}\\
                        \end{array} \right.    
                \end{array}$} & \hspace{-1cm}  \\[30ex]
                \hline\\[-1 ex]

        \footnotesize{\emph{16}} \hspace{-1cm}   & \footnotesize{$\begin{array}{ll}
            &\emph{(A0)} \begin{array}{ll}
                    L\cdot(&\hspace{-0.3 cm}-K_{S_{16}})=\dfrac{7n+2}{(28n-6)(63n-17)}, ~~R\cdot(-K_{S_{16}})=\dfrac{14n+4}{(28n-6)(63n-17)}, ~~L\cdot R=\dfrac{3}{63n-17}\\\\
                    &L^2=\dfrac{7}{(28n-6)(63n-17)}-\dfrac{3}{63n-17}, ~~R^2=\dfrac{14}{(28n-6)(63n-17)}-\dfrac{3}{63n-17}\\\\
                   \end{array} \\\\
            &\emph{(A1)} ~~\dfrac{3(7n+2)}{(28n-6)(63n-17)}=H_x\cdot D \geq \mult_{\mathsf{p}}(H_x)\cdot \mult_{\mathsf{p}}(D)>\left\{
                \begin{array}{ll}
                1, ~~\mathrm{if} ~\mathsf{p} ~~\mathrm{is ~~smooth}\\\\   
                \dfrac{2}{28n-6}, ~~\mathrm{if} ~\mathsf{p}=\mathsf{p}_{_y}\\\\
                \dfrac{2}{63n-17}, ~~\mathrm{if} ~\mathsf{p}=\mathsf{p}_{_t}\\\\
                \dfrac{2}{14n-3}, ~~\mathrm{if} ~\mathsf{p}=\mathsf{p}_{_{14n-3}}\\
                \end{array} \right.\\\\
            &\emph{(A2)} ~~\dfrac{7n+2}{(28n-6)(63n-17)}=L\cdot D > \dfrac{1}{63n-17}
            ~~~\left(\mathrm{resp}. ~~\dfrac{14n+4}{(28n-6)(63n-17)}=R\cdot D > \dfrac{1}{63n-17} \right)\\\\\\
            
            &\emph{(A3)} ~~\dfrac{14n+4+\lambda_{R}(84n-32)}{(28n-6)(63n-17)}=R\cdot (D-\lambda_{R}R)= R \cdot\Delta > \left\{ 
                        \begin{array}{ll}
                        1, ~~\mathrm{if} ~\mathsf{p} ~~\mathrm{is ~~smooth}\\\\   
                        \dfrac{1}{14n-3}, ~~\mathrm{if} ~\mathsf{p}=\mathsf{p}_{_{14n-3}}\\
                        \end{array} \right.    
                \end{array}$} &  \hspace{-1cm} \\[30ex]
                \hline\\[-1 ex]

        \footnotesize{\emph{18}} \hspace{-1cm}  & \footnotesize{$\begin{array}{ll}
            &\emph{(A0)} \begin{array}{ll}
                    L\cdot(&\hspace{-0.3 cm}-K_{S_{18}})=\dfrac{7n+3}{(28n-2)(63n-8)}, ~~R\cdot(-K_{S_{18}})=\dfrac{14n+6}{(28n-2)(63n-8)}, ~~L\cdot R=\dfrac{3}{63n-8}\\\\
                    &L^2=\dfrac{7}{(28n-2)(63n-8)}-\dfrac{3}{63n-8}, ~~R^2=\dfrac{14}{(28n-2)(63n-8)}-\dfrac{3}{63n-8}\\\\
                   \end{array} \\\\
            &\emph{(A1)} ~~\dfrac{3(7n+3)}{(28n-2)(63n-8)}=H_x\cdot D \geq \mult_{\mathsf{p}}(H_x)\cdot \mult_{\mathsf{p}}(D)>\left\{
                \begin{array}{ll}
                1, ~~\mathrm{if} ~\mathsf{p} ~~\mathrm{is ~~smooth}\\\\   
                \dfrac{2}{28n-2}, ~~\mathrm{if} ~\mathsf{p}=\mathsf{p}_{_y}\\\\
                \dfrac{2}{63n-8}, ~~\mathrm{if} ~\mathsf{p}=\mathsf{p}_{_t}\\\\
                \dfrac{2}{14n-1}, ~~\mathrm{if} ~\mathsf{p}=\mathsf{p}_{_{14n-1}}\\
                \end{array} \right.\\\\
            &\emph{(A2)} ~~\dfrac{7n+3}{(28n-2)(63n-8)}=L\cdot D > \dfrac{1}{63n-8}
            ~~~\left(\mathrm{resp}. ~~\dfrac{14n+6}{(28n-2)(63n-8)}=R\cdot D > \dfrac{1}{63n-8} \right)\\\\\\
            
            &\emph{(A3)} ~~\dfrac{14n+6+\lambda_{R}(84n-20)}{(28n-2)(63n-8)}=R\cdot (D-\lambda_{R}R)= R \cdot\Delta > \left\{ 
                        \begin{array}{ll}
                        1, ~~\mathrm{if} ~\mathsf{p} ~~\mathrm{is ~~smooth}\\\\   
                        \dfrac{1}{14n-1}, ~~\mathrm{if} ~\mathsf{p}=\mathsf{p}_{_{14n-1}}\\
                        \end{array} \right.    
                \end{array}$} &  \hspace{-1cm} \\[30ex]
                \hline\\[-1 ex]

        \footnotesize{\emph{20}} \hspace{-1cm}  & \footnotesize{$\begin{array}{ll}
            &\emph{(A0)} \begin{array}{ll}
                    L\cdot(&\hspace{-0.3 cm}-K_{S_{20}})=\dfrac{7n+4}{(28n+2)(63n+1)}, ~~R\cdot(-K_{S_{20}})=\dfrac{14n+8}{(28n+2)(63n+1)}, ~~L\cdot R=\dfrac{3}{63n+1}\\\\
                    &L^2=\dfrac{7}{(28n+2)(63n+1)}-\dfrac{3}{63n+1}, ~~R^2=\dfrac{14}{(28n+2)(63n+1)}-\dfrac{3}{63n+1}\\\\
                   \end{array} \\\\
            &\emph{(A1)} ~~\dfrac{3(7n+4)}{(28n+2)(63n+1)}=H_x\cdot D \geq \mult_{\mathsf{p}}(H_x)\cdot \mult_{\mathsf{p}}(D)>\left\{
                \begin{array}{ll}
                1, ~~\mathrm{if} ~\mathsf{p} ~~\mathrm{is ~~smooth}\\\\   
                \dfrac{2}{28n+2}, ~~\mathrm{if} ~\mathsf{p}=\mathsf{p}_{_y}\\\\
                \dfrac{2}{63n+1}, ~~\mathrm{if} ~\mathsf{p}=\mathsf{p}_{_t}\\\\
                \dfrac{2}{14n+1}, ~~\mathrm{if} ~\mathsf{p}=\mathsf{p}_{_{14n+1}}\\
                \end{array} \right.\\\\
            &\emph{(A2)} ~~\dfrac{7n+4}{(28n+2)(63n+1)}=L\cdot D > \dfrac{1}{63n+1}
            ~~~\left(\mathrm{resp}. ~~\dfrac{14n+8}{(28n+2)(63n+1)}=R\cdot D > \dfrac{1}{63n+1} \right)\\\\\\
            
            &\emph{(A3)} ~~\dfrac{14n+8+\lambda_{R}(84n-8)}{(28n+2)(63n+1)}=R\cdot (D-\lambda_{R}R)= R \cdot\Delta > \left\{ 
                        \begin{array}{ll}
                        1, ~~\mathrm{if} ~\mathsf{p} ~~\mathrm{is ~~smooth}\\\\   
                        \dfrac{1}{14n+1}, ~~\mathrm{if} ~\mathsf{p}=\mathsf{p}_{_{14n+1}}\\
                        \end{array} \right.    
                \end{array}$} &  \hspace{-1cm} \\[30ex]
                \hline\\[-1 ex]

        \footnotesize{\emph{22}} \hspace{-1cm}  & \footnotesize{$\begin{array}{ll}
            &\emph{(A0)} \begin{array}{ll}
                    L\cdot(&\hspace{-0.3 cm}-K_{S_{22}})=\dfrac{7n+5}{(28n+6)(63n+10)}, ~~R\cdot(-K_{S_{22}})=\dfrac{14n+10}{(28n+6)(63n+10)}, ~~L\cdot R=\dfrac{3}{63n+10}\\\\
                    &L^2=\dfrac{7}{(28n+6)(63n+10)}-\dfrac{3}{63n+10}, ~~R^2=\dfrac{14}{(28n+6)(63n+10)}-\dfrac{3}{63n+10}\\\\
                   \end{array} \\\\
            &\emph{(A1)} ~~\dfrac{3(7n+5)}{(28n+6)(63n+10)}=H_x\cdot D \geq \mult_{\mathsf{p}}(H_x)\cdot \mult_{\mathsf{p}}(D)>\left\{
                \begin{array}{ll}
                1, ~~\mathrm{if} ~\mathsf{p} ~~\mathrm{is ~~smooth}\\\\   
                \dfrac{2}{28n+6}, ~~\mathrm{if} ~\mathsf{p}=\mathsf{p}_{_y}\\\\
                \dfrac{2}{63n+10}, ~~\mathrm{if} ~\mathsf{p}=\mathsf{p}_{_t}\\\\
                \dfrac{2}{14n+3}, ~~\mathrm{if} ~\mathsf{p}=\mathsf{p}_{_{14n+3}}\\
                \end{array} \right.\\\\
            &\emph{(A2)} ~~\dfrac{7n+5}{(28n+6)(63n+10)}=L\cdot D > \dfrac{1}{63n+10}
            ~~~\left(\mathrm{resp}. ~~\dfrac{14n+10}{(28n+6)(63n+10)}=R\cdot D > \dfrac{1}{63n+10} \right)\\\\\\
            
            &\emph{(A3)} ~~\dfrac{14n+10+\lambda_{R}(84n+4)}{(28n+6)(63n+10)}=R\cdot (D-\lambda_{R}R)= R \cdot\Delta > \left\{ 
                        \begin{array}{ll}
                        1, ~~\mathrm{if} ~\mathsf{p} ~~\mathrm{is ~~smooth}\\\\   
                        \dfrac{1}{14n+3}, ~~\mathrm{if} ~\mathsf{p}=\mathsf{p}_{_{14n+3}}\\
                        \end{array} \right.    
                \end{array}$} &\hspace{-1cm}  \\[30ex]
        \hline\hline
        \end{longtable}
        \end{center}

        \newpage

        Each absurd inequality \emph{(A4)} where $H_x$ is reducible as in Theorem \ref*{main1} is presented in Table 7.
    \begin{center}
        \begin{longtable}{ccccc}
            \caption{Absurd Inequalities for Complete Intersections}\\
        \hline\\[-1.5 ex]
        \footnotesize{No.} \hspace{0.4cm} &  & \footnotesize{Defining Equations and Absurd Inequalities } & & \hspace{0.4cm} \\ [1 ex]
        \hline\hline\\[-1.5 ex]
    
        \footnotesize{\emph{7}}\hspace{0.4cm} & & \footnotesize{$\begin{array}{ll}
            &\mathbb{P}(1,8n-2,12n-3,18n-5,24n-7), ~~I=2n; ~~\left\{ \begin{array}{ll}
                        wx-z^2-y^3=0\\\\   
                        wz+t^2+x^{36n-10}=0\\
                        \end{array} \right. \\\\
            & \emph{(A4)} ~~\dfrac{4n}{(4n-1)(24n-7)}=R'\cdot D' \geq \mult_{\mathsf{q}}(R')\cdot \mult_{\mathsf{q}}(D')> \dfrac{3}{24n-7}\\
                \end{array}$} & &  \hspace{0.4cm}\\[7ex]
                \hline\\[-1 ex]

        \footnotesize{\emph{11}}\hspace{0.4cm} &  & \footnotesize{$\begin{array}{ll}
            &\mathbb{P}(7,28n-18,42n-27,63n-44,84n-61), ~~I=7n-1; 
            ~~\left\{ \begin{array}{ll}
                        wx-z^2-y^3=0\\\\   
                        wz+t^2+x^{14n-10}y=0\\
                        \end{array} \right. \\\\
            & \emph{(A4)} ~~\dfrac{14n-2}{(14n-9)(84n-61)}=R'\cdot D' \geq \mult_{\mathsf{q}}(R')\cdot \mult_{\mathsf{q}}(D')> \dfrac{3}{84n-61}\\
                \end{array}$} &  & \hspace{0.4cm}\\[7ex]
                \hline\\[-1 ex]

        \footnotesize{\emph{14}}\hspace{0.4cm}  & & \footnotesize{$\begin{array}{ll}
            &\mathbb{P}(7,28n-10,42n-15,63n-26,84n-37), ~~I=7n+1; 
            ~~\left\{ \begin{array}{ll}
                        wx-z^2-y^3=0\\\\   
                        wz+t^2+x^{14n-6}y=0\\
                        \end{array} \right. \\\\
            & \emph{(A4)} ~~\dfrac{14n+2}{(14n-5)(84n-37)}=R'\cdot D' \geq \mult_{\mathsf{q}}(R')\cdot \mult_{\mathsf{q}}(D')> \dfrac{3}{84n-37}\\
                \end{array}$} & & \hspace{0.4cm} \\[7ex]
                \hline\\[-1 ex]
    
        \footnotesize{\emph{16}}\hspace{0.4cm} &  & \footnotesize{$\begin{array}{ll}
            &\mathbb{P}(7,28n-6,42n-9,63n-17,84n-25), ~~I=7n+2; 
            ~~\left\{ \begin{array}{ll}
                        wx-z^2-y^3=0\\\\   
                        wz+t^2+x^{14n-4}y=0\\
                        \end{array} \right. \\\\
            & \emph{(A4)} ~~\dfrac{14n+4}{(14n-3)(84n-25)}=R'\cdot D' \geq \mult_{\mathsf{q}}(R')\cdot \mult_{\mathsf{q}}(D')> \dfrac{3}{84n-25}\\
                \end{array}$} &  & \hspace{0.4cm}\\[7ex]
                \hline\\[-1 ex]
    
        \footnotesize{\emph{18}}\hspace{0.4cm} &  & \footnotesize{$\begin{array}{ll}
            &\mathbb{P}(7,28n-2,42n-3,63n-8,84n-13), ~~I=7n+3; 
            ~~\left\{ \begin{array}{ll}
                        wx-z^2-y^3=0\\\\   
                        wz+t^2+x^{14n-2}y=0\\
                        \end{array} \right. \\\\
            & \emph{(A4)} ~~\dfrac{14n+6}{(14n-1)(84n-13)}=R'\cdot D' \geq \mult_{\mathsf{q}}(R')\cdot \mult_{\mathsf{q}}(D')> \dfrac{3}{84n-13}\\
                \end{array}$} & & \hspace{0.4cm} \\[7ex]
                \hline\\[-1 ex]
    
        \footnotesize{\emph{20}}\hspace{0.4cm}  & & \footnotesize{$\begin{array}{ll}
            &\mathbb{P}(7,28n+2,42n+3,63n+1,84n-1), ~~I=7n+4; 
            ~~\left\{ \begin{array}{ll}
                        wx-z^2-y^3=0\\\\   
                        wz+t^2+x^{14n}y=0\\
                        \end{array} \right. \\\\
            & \emph{(A4)} ~~\dfrac{14n+8}{(14n+1)(84n-1)}=R'\cdot D' \geq \mult_{\mathsf{q}}(R')\cdot \mult_{\mathsf{q}}(D')> \dfrac{3}{84n-1}\\
                \end{array}$} &  & \hspace{0.4cm}\\[7ex]
                \hline\\[-1 ex]
    
        \footnotesize{\emph{22}}\hspace{0.4cm} &  & \footnotesize{$\begin{array}{ll}
            &\mathbb{P}(7,28n+6,42n+9,63n+10,84n+11), ~~I=7n+5; 
            ~~\left\{ \begin{array}{ll}
                        wx-z^2-y^3=0\\\\   
                        wz+t^2+x^{14n+2}y=0\\
                        \end{array} \right. \\\\
            & \emph{(A4)} ~~\dfrac{14n+10}{(14n+3)(84n+11)}=R'\cdot D' \geq \mult_{\mathsf{q}}(R')\cdot \mult_{\mathsf{q}}(D')> \dfrac{3}{84n+11}\\
                \end{array}$} & & \hspace{0.4cm} \\[7ex]
        \hline\hline
        \end{longtable}
        \end{center}

\begin{ack}
    The authors are deeply grateful to Ivan Cheltsov. His suggestion was the starting point to study this problem. The first author was supported by the National Research Foundation of Korea (NRF-2023R1A2C1003390 and NRF-2022M3C1C8094326). The second author was supported by the National Research Foundation of Korea (NRF-2020R1A2C1A01008018, NRF-2022M3C1C8094326 and NRF-2021R1A6A1A10039823). The third authors was supported by the National Research Foundation of Korea (NRF-2020R1A2C1A01008018 and NRF-2022M3C1C8094326).
\end{ack}

\end{document}